\documentclass[12pt]{amsart}

\usepackage{amssymb}
\usepackage{fancyhdr}
\usepackage{amsmath}
\usepackage{mathrsfs}
\usepackage{amsfonts}
\usepackage[all]{xy} 
\usepackage{xcolor}

\newcommand{\disk}{\ensuremath{\mathbb{D}} } 
\newcommand{\sphere}{\bar{\Bbb{C}}} 
\newcommand{\riem}{\Sigma}  
\renewcommand{\Bbb}[1]{\ensuremath{\mathbb{#1}}}

\newcommand{\Oqc}{\mathcal{O}^{\mathrm{qc}}} 
\newcommand{\qs}{\operatorname{QS}}

\newcommand{\qc}{\operatorname{QC}}

\theoremstyle{plain}
        \newtheorem{theorem}{Theorem}[section]
        
        \newtheorem{proposition}[theorem]{Proposition}
        \newtheorem{corollary}[theorem]{Corollary}

\theoremstyle{definition}
        \newtheorem{definition}[theorem]{Definition}
        \newtheorem{example}{Example}[section]

\theoremstyle{remark}
    \newtheorem{remark}[theorem]{Remark}
    \newtheorem{question}{Problem}[section]

\numberwithin{equation}{section} 
\numberwithin{figure}{section} 

\usepackage{fullpage}

\title{Infinite-dimensional moduli spaces of Riemann surfaces}
\dedicatory{Dedicated to the memory of Alexander Vasil'ev}
\begin{document}
\title[Comparison Moduli spaces of Riemann surfaces]{Comparison moduli spaces of Riemann surfaces}

\author{Eric Schippers}
\address{Eric Schippers \\ Department of Mathematics \\
University of Manitoba\\
Winnipeg, Manitoba \\  R3T 2N2 \\ Canada}
\email{eric\_schippers@umanitoba.ca}
\thanks{}

\author{Wolfgang Staubach}
\address{Wolfgang Staubach\\ Department of Mathematics\\
Uppsala University\\
Box 480\\ 751 06 Uppsala\\ Sweden}
\email{wulf@math.uu.se}


\date{\today}

\begin{abstract}
 We define a kind of moduli space of nested surfaces and mappings, which we call a comparison
 moduli space.  We review examples of such spaces in geometric function theory and modern Teichm\"uller
 theory, and illustrate how a wide range of phenomena in complex analysis are captured by this
 notion of moduli space.  The paper includes a list of open problems in classical and modern function theory and Teichm\"uller theory ranging from  general
 theoretical questions to specific technical problems.
\end{abstract}

\thanks{Eric Schippers and Wolfgang Staubach are grateful for the financial support from the Wenner-Gren Foundations. Eric Schippers is also partially supported by the National Sciences and Engineering Research Council of Canada. }

\maketitle

\begin{section}{Introduction}
 Theorems in complex analysis should be conformally invariant.  In any branch
 of mathematics, one does not distinguish objects with respect to some equivalence; in
 complex analysis, this is conformal equivalence.
  Thus, if one alters all objects
 in the hypothesis of a theorem by applying a biholomorphism, then the theorem should apply to the new situation.

 That is not to say that theorems which are not manifestly conformally invariant are not interesting.
 Indeed any pair of points on two distinct Riemann surfaces are contained in locally biholomorphic neighbourhoods
 by the Riemann mapping theorem, so local phenomena need not ever refer to conformal invariance.
 Furthermore we have powerful uniformization theorems at our disposal which reduce problems involving
 Riemann surfaces or mappings between them to
 canonical cases.   This accounts for the fact that many deep results in complex analysis
 need not mention conformal invariance.  However, it is still profitable to attempt to formulate
 - or reformulate - complex analytic results invariantly.\\

 It is not always obvious how to do this.  In the case of results about compact Riemann surfaces, it is
 straightforward.  These Riemann surfaces are characterized by a finite set of conformal invariants,
 or by a finite-dimensional function space which is easily defined conformally invariantly, such
 as the vector space of holomorphic abelian differentials.
 However, things change radically when the Riemann surfaces have boundary curves.  For example, the function spaces
 become infinite-dimensional, boundary values might enter their definition, and furthermore one must become concerned with boundary regularity of the biholomorphism.
 Thus many important complex analytic phenomena do not easily fit into the framework which is so successful
 for compact surfaces.

 An example will illustrate the point.  The universal Teichm\"uller space can be identified
 with the set of suitably normalized domains in the sphere bounded by quasicircles, and can
 be thought of as the set of deformations of the complex structure of the unit disk.  Consider the following reasonable objection to this definition, encountered occasionally by the authors:
 by the Riemann mapping theorem, these domains (or deformations) are all conformally equivalent.
 So why would one care about such a space?  In other words, by conformal invariance, the universal
 Teichm\"uller space is of little interest.  However, if one takes this objection too seriously, one
 will never find out that the universal Teichm\"uller space contains all Teichm\"uller spaces, or
 learn the beautiful geometric and analytic theory of quasidisks and quasiconformal mappings.
 Thus this naive approach to conformal invariance is not appropriate for quasiconformal
 Teichm\"uller theory.

 How then do we introduce conformal invariance into this picture?
 This is done by taking into account the fact that the universal Teichm\"uller space involves two Riemann surfaces:
 the quasidisk and the sphere itself.  Thus it is a moduli space of embeddings of disks into the
 sphere. \\

 In this paper, we outline a general notion of moduli spaces of pairs of surfaces,
 or the closely related moduli spaces of mappings.   We call them
 comparison moduli spaces.
This notion of moduli spaces and the corresponding
 notion of conformal invariance capture a wide range of complex analytic phenomena.
 We show how comparison moduli spaces are present in classical geometric function theory, and along the
 way draw attention to some deep ideas of Z. Nehari and M. Schiffer which in our opinion have received
 too little attention.  We also review some modern manifestations of the idea in quasiconformal
 Teichm\"uller theory and conformal field theory, with special attention to the
 so-called Weil-Petersson class Teichm\"uller space.   The paper also contains a list of open problems.

  In attempting to convey a unifying viewpoint, we
  have chosen to paint a broad rather than a specific picture.
  Thus, our choice of topics is selective; especially so when we are
  dealing with older topics.  We have not attempted to make a survey of current or past research
  on each topic, and have confined ourselves to providing a few key references which
  should be viewed as entry points to the literature.   However, since the  Weil-Petersson class Teichm\"uller
  theory of non-compact surfaces is rather new, we made an exception for this case and gave a more comprehensive literature review.  \\

The organization of this paper is as follows. In Section \ref{se:comparison_moduli_spaces}, we define and give examples of comparison moduli spaces and the corresponding notion of conformal invariance.  In Section \ref{se:geometric_function_theory}
 we review examples of such moduli spaces in geometric function theory, in the work of
 Nehari and Schiffer.  In Section \ref{se:Teich_theory}
 we give a brief overview of quasiconformal Teichm\"uller theory, and the rigged moduli space arising in
 conformal field theory.  We describe the correspondence between these two, which also demonstrates that
 Teichm\"uller space is, up to a discontinuous group action, a comparison moduli space.  Finally, in
 Section \ref{se:WeilPetersson}, we give an overview of new developments
 in Weil-Petersson Teichm\"uller theory. We also
 briefly list some of the applications of this new field.

\end{section}
\textbf{Acknowledgements.} The authors are grateful to David Radnell for a fruitful collaboration and valuable discussions through the years, and for his comments and suggestions concerning the initial draft of the manuscript.
\begin{section}{Comparison moduli spaces}\label{se:comparison_moduli_spaces}

 By a comparison moduli space, we mean a moduli space of one of two types, which
 we now describe.  Our purpose is to illustrate the similarity of a wide range of
 ideas, and how the notion of comparison moduli space is useful in capturing complex
   analytic phenomena.  Our purpose is {\it not} to give a Bourbaki styled axiomatization.
 Although the definitions of the comparison moduli spaces
 can easily be made precise, this introduces
 unnecessary abstractions.  Thus, rather than taking this approach, we give them informal but general definitions. Precise definitions are reserved for the many examples throughout the paper.

\begin{subsection}{Moduli space of nested surfaces}
  Let $R_1$ and $R_2$ be Riemann surfaces such that $R_2 \subseteq R_1$.
  We will also usually equip the nested surfaces with extra data $a,b,c,\ldots$.
  We assume that for any biholomorphism $g:S_1 \rightarrow R_1$, there is an operation
  $g^*$ on the data denoted $a \mapsto g^*(a)$ which takes data on
  $(R_1,R_2)$ to data on $(S_1,S_2)$.
  For example, $a$ might be a specified point in $S_1$ and $g^*(a) = g^{-1}(a)$,
  while $b$ might be a harmonic function on $S_1 \backslash \{ a \}$ and
  $g^* b = b \circ g$.  This operation should be invertible and respect composition:
  that is, $(g_1 \circ g_2)^*(a)= g_2^*(g_1^*(a))$ and $(g^{-1})^* =  (g^*)^{-1}$.
  We call this operation the {\it pull-back} under $g$, and call the set $\{(R_1,R_2,a,b,c,\ldots) \}$
  the {\it configuration space} $\mathfrak{C}$.

  Two elements of the configuration space are equivalent
  \[  (R_1,R_2,a,b,c,\ldots)  \sim (S_1,S_2,A,B,C,\ldots) \]
  if and only if there is a biholomorphism $g:S_1 \rightarrow R_1$
  such that
  \[  g(S_2)=R_2, \ \ g^*(a)=A, \ \ g^*(b)=B, \ \ g^*(c)=C,\ldots  \ .  \]
   The {\it moduli space of nested surfaces} $\mathfrak{M}$ is the set of equivalence classes
  \[  \mathfrak{M} = \mathfrak{C} / \sim. \]
  A {\it conformal invariant} is a function on the configuration space
  $I:\mathfrak{C} \rightarrow \mathbb{C}$ (or into $\mathbb{R}$)
  such that for any biholomorphic map $g:S_1 \rightarrow R_1$ such that
  $g(S_2)=R_2$ one has
  \begin{equation}  \label{eq:conf_inv_nested_surfaces}
   I(R_1,R_2,a,b,c,\ldots) = I(S_1,S_2,g^*(a),g^*(b),g^*(c),\ldots).
  \end{equation}
    Equivalently, it is a function on $\mathfrak{M}$.
  \begin{remark} \label{re:extension_to_differentials}  The definition \ref{eq:conf_inv_nested_surfaces} can be generalized
   easily to allow for the possibility that $I$ is a differential of order $n$
   on $R_2$ (that is, given in local coordinates by $h(z)dz^n$).
  \end{remark}

  Finally, one might place restrictions on the Riemann surfaces (e.g.
  bounded by quasicircles, fixed genus) or the data (e.g. harmonic
  functions of finite Dirichlet energy, singularity of specified form
  at a point).  It is required that these conditions are conformally invariant.

  \begin{example}  \label{ex:hyperbolic_surfaces_lambda_ratio}
   Let $\mathfrak{C} =\{(R_1,R_2,z)\}$ such that
   $R_1$ and $R_2$ are hyperbolic Riemann surfaces, $R_2 \subseteq R_1$ and $z \in R_2$.
   Let $\lambda_i(z)^2|dz|^2$ denote the hyperbolic metric on $R_1$ in any fixed
   local coordinate (with curvature normalized
   to $-1$ say).  For a biholomorphism $g:S_1 \rightarrow R_1$ define $g^*(z)=g^{-1}(z)$.
   Define
   \[  I(R_1,R_2,z) = \lambda_1(z)/\lambda_2(z).  \]
   This is independent of the choice of coordinate and conformally invariant.

   Another example of a conformal invariant is the following \cite{Schippers_confinv_analyse}.
   Defining $\Gamma_i = 2 \frac{\partial}{\partial z} \log{\lambda_i}$ in some
   local coordinate, the quantity
   \[  I_2(R_1,R_2,z) = \lambda_2^{-1}(z) | \Gamma_{2}(z) - \Gamma_{1}(z) |  \]
   is independent of choice of coordinate and conformally invariant.
  \end{example}
  \begin{example}  Let $\mathfrak{C} =\{(R_1,R_2)\}$ where
   $R_1$ is a simply connected Riemann surface which is biholomorphic
   to the disk and $R_2 \subseteq R_1$ is bounded by a quasicircle in $R_1$.
   The condition that $R_2$ is bounded by a quasicircle in $R_1$ is conformally
   invariant.
  \end{example}
  \begin{example}  Let $\mathfrak{C} = \{ (R_1,R_2,h) \}$ where
  $R_1$ and $R_2$ are doubly connected Riemann surfaces (i.e.
   biholomorphic to $r<|z|<R$ for $0 < r < R < \infty$), $R_2 \subseteq R_1$
   and $h$ is a harmonic function on $R_1$ of finite Dirichlet energy:
   \[  \iint_{R_1} dh \wedge \overline{dh} <\infty.  \]
   These conditions are conformally invariant. For a biholomorphism $g:S_1
   \rightarrow R_1$ define $g^*h = h \circ g$.  The quantity
   \[  I(R_1,R_2,h) = \iint_{R_1 \backslash R_2} dh \wedge \overline{dh} \]
   is a conformal invariant.
  \end{example}
  \begin{example} \label{ex:univ_Teich_nested} Let $\mathfrak{C} =\{(R_1,R_2,a,b,c)\}$
   where $R_1$ is a Riemann surface biholomorphic to the Riemann sphere,
   $R_2 \subseteq R_1$ is bordered by a quasicircle in $R_1$, and $a,b,c$ are
   points on $\partial R_2$.  For a biholomorphism $g:S_1 \rightarrow R_1$
   we define $g^*(a)=g^{-1}(a)$, and similarly for $b$ and $c$.  The moduli space
   $\mathfrak{M}$ is then the universal Teichm\"uller space.  This will be discussed in Section \ref{se:qc_Teich_primer} ahead.
  \end{example}
\end{subsection}
\begin{subsection}{Moduli spaces of mappings}  \label{se:moduli_spaces_mappings}
 We consider also moduli spaces of mappings between Riemann surfaces.
 That is, we define a configuration space $\mathfrak{C} = \{ (R_0,f,R,a,b,c,\ldots) \}$
 where $R_0$ is a fixed Riemann surface, $R$ is a variable Riemann surface, $f:S \rightarrow R$ is
 an injective holomorphic map, and $a,b,c,\ldots$ are data as above.
 Pull-backs under biholomorphisms are defined as above. We define
 an equivalence relation  $(R_0,f,R,a,b,c,\ldots)\sim (R_0,F,S,A,B,C,\ldots)$
 if and only if there is a biholomorphism $g:S \rightarrow R$ such that
 $f = g \circ F$, $g^*(a)=A$, $g^*(b)=B$, $g^*(c)=C, \ldots$  .
 The moduli space of mappings is
 \[  \mathfrak{M} = \mathfrak{C} /\sim.  \]
 A conformal invariant is a quantity $I$ on $\mathfrak{C}$ such that for
 any biholomorphism $g:S \rightarrow R$
 \begin{equation} \label{eq:map_invariance}
  I(R_0,g \circ F,R,a,b,c,\ldots) = I(R_0,F,S,g^*(a),g^*(b),g^*(c),\ldots).
 \end{equation}
 Equivalently, it is a function on $\mathfrak{M}$.

 We give further examples.  Denote the unit disk by $\disk = \{z \,:\,|z|<1 \}$, the complex plane by
 $\mathbb{C}$, and the Riemann sphere by $\sphere$.
 \begin{example} Let $\mathfrak{C}=\{(\disk,f,R,z) \}$
  where $R$ is a simply connected Riemann surface conformally equivalent to $\disk$,
  $f$ is an injective holomorphic mapping such that $f(0)=z$.
  Then $\mathfrak{M}$ is in one-to-one correspondence with the set of bounded
  univalent functions $f:\disk \rightarrow \disk$ such that $f(0)=0$.

  Note that if we set $R_1 = R$ and $R_2 = f(\disk)$ then
  $\mathfrak{C}$ and $\mathfrak{M}$
  corresponds to the comparison moduli spaces of Example \ref{ex:hyperbolic_surfaces_lambda_ratio},
  modulo an action by $\mathbb{S}^1$ (since the domain $f(\disk)$ does not uniquely determine $f$).
  Setting $z=0$ and $R = \disk$ we can compute the conformal invariant $I_2$  explicitly
  in terms of $f$, namely $I_2(\disk,f,R,z)=\left|f''(0)/f'(0)^2\right|$.  Since this is invariant
  under the $\mathbb{S}^1$-action $f(z) \mapsto e^{-i\theta}f(e^{i\theta}z)$, the invariant
  $I_2$ is well-defined on the configuration $\mathfrak{C}$ considered here.
  The general expression in terms of $z$ is more complicated.
 \end{example}
 \begin{remark}  The notion of comparison moduli spaces is required in order to formulate statements regarding higher-order derivatives of maps conformally invariantly \cite{Schippers_quadratic_expo}.
 \end{remark}
 \begin{example}  \label{ex:univ_Teich_mapping}
  Let $\mathfrak{C} = \{(\disk,f,R) \}$ where $R$ is conformally equivalent to the
  sphere $\sphere$, $a,b,c$ are points in $R$, and  $f:\disk \rightarrow R$ is an injective holomorphic map with a quasiconformal
  extension to $R$.    Then $\mathfrak{M}$ is easily seen to be in one-to-one correspondence with Example \ref{ex:univ_Teich_nested}, and thus is the universal Teichm\"uller space.  Other choices of normalization on $f$ can be imposed.
 \end{example}
 \begin{example}  \label{ex:rigged_moduli_space_analytic}
  Let $\mathfrak{C} = \{(\disk,f_1,\ldots,f_n,R,p_1,\ldots,p_n) \}$ where $R$
  is a compact Riemann surface of genus $g$, $p_1,\ldots,p_n \in R$, $f_i:\disk \rightarrow R$
  such that $f_i$ are one-to-one holomorphic maps on the closure of the disk
  $\text{cl} \disk$, the closures of the images of $\disk$ do not overlap,
  and $f_i(0)=p_i$.  Then $\mathfrak{M}$ is a moduli space in conformal field theory
  due to D. Friedan and S. Shenker, G. Segal, and C. Vafa, see Y-Z. Huang \cite{Huang}
  for references.
  The case where analyticity on $\text{cl}\mathbb{D}$ is weakened to quasiconformal extendibility will be
  considered at length in this paper.
 \end{example}

 As one can also see from these examples, given a moduli space of nested surfaces,
 there is a corresponding moduli space of mappings, if one adds the correct data.
\end{subsection}
\end{section}
\begin{section}{Some examples in geometric function theory} \label{se:geometric_function_theory}
\begin{subsection}{Nehari monotonicity theorems and generalizations}
 Classically, Riemann surfaces are characterized by their spaces of
 functions or differentials with specified singularities.  The existence of functions or differentials with
 specified singularities is closely related to the Dirichlet problem
 on Riemann surfaces. In the case of compact surfaces, algebraic geometric
 techniques nearly (but not quite) eliminate the need to consider the Dirichlet problem.
 On the other hand, for nested surfaces with boundaries, there is not
 sufficient rigidity for these algebraic geometric techniques to determine the families of differentials,
 and the Dirichlet problem again takes centre stage.

 Z. Nehari \cite{Nehari_some_inequalities} defined a class of functionals on subdomains of Riemann surfaces,
 obtained from the Dirichlet integrals of harmonic functions of specified
 singularities.    It is closely related to the so-called contour integral and
 area techniques in function theory.  We will not review these connections
 since they can be found in other sources.

 Nehari considered the comparison moduli space $\mathfrak{C}$ (he did not
 use the term), which in our terminology is given by
 $\mathfrak{C} = \{ (R_1,R_2,z_1,z_2,\ldots,z_n,h) \}$ where $R_1$ and $R_2$ are Riemann surfaces,
 $R_2 \subseteq R_1$ is bounded by analytic curves in the interior of $R_1$, $z_1,\ldots,z_n
 \in R_2$ and $h$ is a harmonic function on $R_1$ with singularities at $z_i$.  He did
 not precisely specify the nature of the singularities of $h$, but in all examples they were such that $h$
 was locally the real part of a meromorphic function with poles at $z_i$.  If
 we define pull-back on points and the harmonic function in the obvious way, these conditions
 are conformally invariant and thus we obtain the moduli space $\mathfrak{M}$.\\

Below, we shall give a collection of monotonic functionals on $\mathfrak{M}$ defined by Nehari. We will also shorten the notation of the configuration space to $(R_1,R_2,h)$ (since
 the information of the location of the singularities are encoded in the domain
 of the harmonic function $h$).
 \begin{definition}[Nehari functional]
  Let $R$ be a Riemann surface, and let $D$ be a
  subdomain of $R$ bounded by finitely many analytic curves.
  Let $h$ be a harmonic function with a finite number of isolated singularities
  all of which are in $D$.
  Let $q$ be the unique function on $D$ which is constant on
  each component of $\partial D$, such that $q + h$ is harmonic on $D$,
  and such that for any closed contour $\Gamma$ in $D$,
  \[  \int_{\Gamma} \frac{\partial(q + h)}{\partial n} ds =0. \]
 Here $n$ denotes the unit outward normal and $ds$ denotes the unit arc length.
  We define the ``Nehari functional'' as
  \[  N(R,D,h) = \int_{\partial D} h \frac{\partial q}{\partial n} ds.  \]
 \end{definition}
 \begin{remark}  Although the outward normal is coordinate dependent, the expression
  $h \frac{\partial q}{\partial n} ds$ is not.
 \end{remark}
 It is immediately evident that the Nehari functional is conformally invariant
 in the sense that if $g:R \rightarrow R'$ is a biholomorphism of Riemann
 surfaces then $N(R,D,h)=N(R',g(D),h \circ g^{-1})$.
 \begin{theorem}[Nehari \cite{Nehari_some_inequalities}]
  $N$ is a monotonic functional in the sense that whenever $D_1 \subseteq D_2$
 it holds that $N(R,D_2 ,h) \leq N(R,D_1 ,h)$.
  \end{theorem}
   This follows from the Dirichlet principle.  It can be generalized to collections of non-overlapping domains.

 Nehari showed that this concept unifies a large number of results in function
  theory, including the Grunsky inequalities, many coefficient estimates for
  univalent functions, and estimates on capacity.  Diverse applications
  from the point of view of capacitance
  can also be found in the monograph of V. Dubinin \cite{Dubinin_book}.

  We restrict now to the case that both domains are simply connected
  and hyperbolic. In \cite{Schippers_confinv_analyse,Schippers_confinv_Jyv}, E. Schippers obtained many estimates on conformal invariants involving bounded univalent functions from Nehari's method.
  The order of the estimates are determined by the order of the
  singularity.
   By Teichm\"uller's principle \cite{Jenkins_book},
   extremal problems are in general
  associated with quadratic differentials (this principle was enunciated by J. A. Jenkins who attributes it to O. Teichm\"uller).  If the estimate involves $f^{(n)}(z)$
  at a fixed point $z$, then the associated quadratic differential has a pole
  of order $n+1$ at $z$.

  Thus, one is led to\begin{enumerate}
    \item formulate Nehari's functional in terms of quadratic
     differentials.
  \item consider quadratic differentials with poles of arbitrary order.
  \end{enumerate}

  In constructing this reformulation in terms of quadratic differentials, we restrict to the case that the outside domain $R$ is conformally equivalent to a disk.  Since we are dealing with domains with a global coordinate,
  we will use the notation $Q(z)dz^2$ for quadratic differentials.
  In \cite{Schippers_Israel_Journal} E. Schippers showed that, by altering $h$, every Nehari functional is up to a constant
  equal to a Nehari functional
  such that the boundary of $R$ is a trajectory of
  the quadratic differential $Q(z)dz^2 = h'(z)^2dz^2$.
    Furthermore, this particular choice
  of $h$ is precisely the one so that $N(R,R,h)=0$.

  Of course, this strongly suggests
  that the Nehari functional can be generalized to arbitrary quadratic differentials, not just those of the form $h'(z)^2 dz^2$.
  This can indeed be done by passing to a double cover $\pi:\tilde{R} \rightarrow R$ branched at odd-order poles and zeros \cite{Schippers_Israel_Journal,Schippers_quadratic_expo}.   For any
  quadratic differential $Q(z)dz^2$ such that $\partial R$ is a trajectory, define
  \begin{equation} \label{eq:contour_integral}
   m(R,D,Q(z)dz^2) =  \int_{\partial \tilde{D}} h \frac{\partial q}{\partial n} ds
  \end{equation}
  where $\tilde{D} = \pi^{-1}(D)$, $h =\text{Re} \int \sqrt{Q(\pi(z)) \pi'(z)^2} \,dz$, and $q$ is the unique harmonic function on $\tilde{D}$ such that $q = h$ on $\partial \tilde{D}$.  We also assume that the boundary of $D$
  is sufficiently regular, but the functional extends to the case that $D$ is
  bounded by a Jordan curve.
  Some effort is required to make this definition rigorous; one must deal
  with branch points, show that this integral is single-valued and independent
  of choice of branch of the square root, etc. Once all these issues are settled, one obtains the following result:
  \begin{theorem}[\cite{Schippers_Israel_Journal}]  \label{th:generalized_Nehari}
   Let $R$ and $D$ be conformally equivalent to the unit disk, such that
   $D \subseteq R$, let
   $Q(z)dz^2$ be a quadratic differential such that $\partial R$ is a
   trajectory, all of whose poles are contained in $D$.  Then the following statements hold:\\

   \begin{enumerate}
    \item $m(R,D,Q(z)dz^2)$ is conformally invariant in the sense that
     if $g:R \rightarrow R'$ is a conformal bijection and $g(D)=D'$ then
     $m(R',D',Q(z)dz^2) = m(R,D,Q(g(z))g'(z)^2dz^2)$;\\

    \item if $D_1 \subseteq D_2$ then $m(R,D_1,Q(z)dz^2) \leq m(R,D_2,Q(z)dz^2)$;\\

    \item for any $D$, $m(R,D,Q(z)dz^2) \leq 0$;\\

    \item if $m(R,D,Q(z)dz^2) = 0$ then $D$ is $R$ minus arcs of trajectories
    of $Q(z)dz^2$;\\

    \item if $Q(z)dz^2 = h'(z)^2 dz^2$ for some harmonic function $h$
    with singularities, then \linebreak $m(R,D,Q(z)dz^2) = (-1/2) N(R,D,h)$.
   \end{enumerate}
  \end{theorem}
  \begin{remark}  These invariants can be extended to domains slit by analytic arcs.
  \end{remark}

  By choosing $R= \disk$ and $D = f(\disk)$ for a univalent function $f$,
  many inequalities for bounded univalent functions follow \cite{Schippers_Israel_Journal,Schippers_quadratic_expo}.  The functional
  associated with $Q(z)dz^2$ is in accordance with Teichm\"uller's principle.  That is,
  if the quadratic differential has a pole of order $n+1$ at $z$, the $n$th derivative
  of the mapping function $f^{(n)}(z)$ arises in the functional.

  We now pose the following problem:
  \begin{question}  \label{qu:Nehari_to_qc} Generalize this to the case that the
   inner domain is bordered by quasicircles.
  \end{question}
  The geometric motivation for this problem will be discussed ahead; see
  Problems \ref{qu:extend_to_rigged}, \ref{qu:extend_to_TA}, and \ref{qu:extend_to_Teich}.

  The next problems are much harder, and are at the root
  of many problems in geometric function theory.  In the remainder of
  the section, we restrict to the case of two
  nested simply connected hyperbolic Riemann surfaces (although it is clear that
  the problems have natural generalizations).  By conformal invariance
  we may choose $R = \disk$.  We have that any conformal
  map from $\disk$ into $\disk$ which is admissible for a quadratic differential
  $Q(z)dz^2$ attains the upper bound of the functional $m(\disk,f(\disk),Q(z)dz^2)$.

  Now suppose that we are given a functional in advance, which is not obviously of the form $m(\disk,f(\disk),Q(z)dz^2)$ - for example, a coefficient functional for bounded univalent functions. How does one choose a
  quadratic differential so that $m(\disk,f(\disk),Q(z)dz^2)$ is the desired functional?
  Thus we have two problems.
  \begin{question}  Let $\Phi_{Q(z)dz^2}(f) = m(\disk,f(\disk),Q(z)dz^2)$.
   Describe the map $Q(z)dz^2 \mapsto \Phi_{Q(z)dz^2}$ algebraically.
  \end{question}

  Conversely we ask, in the same spirit,
  \begin{question}
   Describe the set of functionals arising from quadratic differentials (say
   with a pole at the origin, and no other poles).  Describe the inverse map
   $\Phi_{Q(z)dz^2} \mapsto Q(z)dz^2$.
  \end{question}

 There are many interpretations of ``algebraically'', of course.
    A satisfactory answer will almost certainly involve the Lie-theoretic properties
     of the classes of (sufficiently regular) bounded and unbounded univalent
  functions; see e.g. A. A. Kirillov and D. V. Yuri'ev \cite{KY2},
  I. Markina, D. Prokhorov and A. Vasil'ev
  \cite{Markina_Prokhorov_Vasilev}, and Schippers \cite{Schippers_power_quadratic}.

  The heart of the problem
  is that the relation between quadratic differentials and functionals has been
  imprecise since its introduction by Schiffer.
  It is a hard problem to construct a method for producing inequalities.   It is
  still harder to obtain an inequality which one decides on at the outset, no matter which of the existing methods is used.
  This problem has never
  been satisfactorily resolved, in spite of the many successes of the theory of extremal problems.

  It has long been known that there is a correspondence between quadratic differentials, boundary points of function spaces, and extremals of functions.  The following references could serve as a starting point:
  P. Duren \cite{Durenbook}, C. Pommerenke \cite{Pommerenkebook}, A. Schaeffer and D. Spencer  \cite{SchaefferSpencer}, and O. Tammi \cite{Tammi_v1,Tammi_v2}.   Although there are many heuristic principles and concrete theorems, the precise relation has never been established.    This may be partly because the proof of the Bieberbach conjecture
  by deBranges has had the unfortunate effect of diverting attention from the problem (unjustly, since
  the Bieberbach conjecture involves only a single point on the boundary).

  A full exposition of these problems
  would take more space than we have here, so we will content ourselves by quoting an elegant result of
  A. Pfluger.  Since this result strongly suggests that much more remains to be discovered, we will follow it with a few related problems.
  \begin{theorem}[\cite{Pfluger2}]  Consider the class of mappings $f(z) = z + a_2 z^2 + a_3 z^3 + \cdots$
   where $f:\disk \rightarrow \mathbb{C}$ is holomorphic and one-to-one.  For every $\lambda \in \mathbb{C} \backslash \{1\}$, there are precisely two mappings $f_\lambda(z)$
   and $-f_\lambda(-z)$ maximizing the functional $f \mapsto \text{Re}(a_3 - \lambda a_2^2)$.  For every $A \in \mathbb{C} \backslash \{0\}$, with the exception of the Koebe function, every function $f$ mapping onto the plane minus trajectories of the quadratic differential
   $(1+Aw)w^{-4} dw^2$ is extremal for precisely one such functional.
  \end{theorem}

  The class of functions is closely related to the universal Teichm\"uller space.  In fact, the set of such functions  which are additionally supposed to be quasiconformally extendible can be identified with the Bers fibre space over the Teichm\"uller space \cite{Takhtajan_Teo_Memoirs}, which is also a comparison moduli space.  In some sense Pfluger's result gives information about the boundary of this moduli space.  Some natural questions follow.

  \begin{question} Does Pfluger's result hold for higher-order functionals?  For example, is there
   an analogous result for functionals of the form $\text{Re}( a_4 + \lambda a_3 a_2 + \mu a_2^3)$
   and quadratic differentials of the form $((1+Aw + Bw^2)/w^5) \,dw^2$?
  \end{question}
  The form of the functional is obtained by demanding that the expression be homogeneous with
  respect to the transformation $f(z) \mapsto e^{-i\theta} f(e^{i\theta} z)$, while the form of
  the quadratic differential can be obtained using Schiffer's variational theorem to
   differentiate the functional \cite[Chapter 7]{Pommerenkebook}.

  Similarly, we can ask whether Pfluger's results hold in some sense for the comparison moduli space
  considered throughout this section:
  \begin{question}  Is there a version of Pfluger's result for bounded one-to-one holomorphic mappings
   $f:\disk \rightarrow \disk$ such that $f(0)=0$?
  \end{question}

  We observe that Pfluger's result was obtained by combining Schiffer variation with J. Jenkins'
  quadratic differential/extremal metric methods.  Possibly the generalization of
  Nehari's results to quadratic differentials might shed further light on these questions.

  Finally, we state some general problems on the relation to Schiffer variation.
  The Schiffer variational method (in its various forms) says that for
  various classes of conformal maps, the extremal
  map of a (reasonably regular) fixed functional maps onto the target
  domain minus trajectories of a quadratic differential \cite{Durenbook,Pommerenkebook}.
  However,
  the quadratic differential is not precisely determined and depends
  to some extent on the unknown extremal function.  In some sense, the Schiffer
  variational method produces the functional derivative of a functional at
  an extremal and relates it to the quadratic differential
  (see \cite[Theorem 7.4]{Pommerenkebook}).

  In the case of the invariants (\ref{eq:contour_integral}), the extremal
  condition suggests that the quadratic differential obtained by Schiffer variation
  equals the original quadratic differential.
  We are thus led to ask:
  \begin{question}  What is the relation between the functional derivative
   of $m(\disk,f(\disk),Q(z)dz^2)$ and the quadratic
   differential $Q(z)dz^2$?
  \end{question}
 By using the Loewner method to obtain functional derivatives, it was shown in \cite{Schippers_Nehari_derivative} that in the case of
  Nehari functionals (that is, for quadratic differentials of the form $h'(z)^2 dz^2$)
  with poles only at the origin, the functional derivative at an extremal is the pull-back
  of the original quadratic differential under the extremal map.  We conjecture that
  this holds in much greater generality, and also ask whether a stronger
  statement can be made.

  \end{subsection}
\begin{subsection}{A few remarks on the extremal metric method}

  J. A. Jenkins' extremal metrics method associates functionals to quadratic differentials \cite{Jenkins_book}.
  This method is based on length-area inequalities of H. Gr\"otzsch \cite{Grotzsch2}
  and Teichm\"uller \cite{Teichmuller_ungleichungen}, which grew into the theory of
  extremal length.  A modern
  exposition can be found in the monograph of A. Vasiliev \cite{Vasilev_monograph}.  According
  to Jenkins \cite{Jenkins_book}, one can interpret these estimates as involving reduced modules;
  see H. G. Schmidt \cite{Schmidt} for an enlightening exposition in special cases.
  However, in practice, this statement is difficult to make precise except
  for quadratic differentials of order two or lower (see e.g. Dubinin
   \cite[Chapter 6]{Dubinin_book}).
  Also, Jenkins' method does not manifestly produce conformally invariant functionals
  for quadratic differentials of order strictly greater than two.   Thus one is
  led to the following problems:
  \begin{question}  Find a systematic definition of conformally invariant
   reduced modules in terms of extremal length of curve families, which
   produces estimates on arbitrary-order derivatives of conformal maps.
  \end{question}
  This is deliberately imprecise, since it has not even been done in the case of bounded univalent maps.
  Here is a more precise formulation in that special case.
  \begin{question}  Let $D_2 \subseteq D_1$ be simply connected domains
   and let $z$ be a point in $D_2$.  Define conformally invariant
   reduced modules of curve families which for $D_1 = \mathbb{D}$
   and $z=0$, explicitly involving the $n$th derivative of the
   conformal map $f:\disk \rightarrow D_1$ taking $0$ to $0$.
   Prove estimates on these modules using extremal length methods of
   Gr\"otzsch/Teichm\"uller/Jenkins.
  \end{question}
  The boundaries of the domains could be chosen as regular as necessary,
  although we note in passing that the quasiconformal invariance of the
  notion of prime ends should allow a quite general formulation.

  These reduced modules should be associated with quadratic differentials
  through their curve families.
  \begin{question}  \label{qu:quad_diff_association}
   Systematically associate a conformally invariant reduced module to
   each quadratic differential (with poles of arbitrary order).
  \end{question}
  Again, this can be formulated for various classes of maps or pairs of domains,
  but it would be of great interest even for bounded univalent maps/simply connected
   domains.

  In some sense Theorem \ref{th:generalized_Nehari} of the first author is an answer
  to Problem \ref{qu:quad_diff_association}, but with Nehari's Dirichlet
  energy technique replacing extremal length.  Now by a theorem of A. Beurling, the
  Dirichlet energy of a mixed Dirichlet/Neumann boundary value problem is the reciprocal of the extremal
  length of a curve family \cite{Ahlfors_conf_inv}.
  This immediately leads to the following problem.

  \begin{question}
   Formulate and prove a generalization of Beurling's theorem
   for harmonic functions with singularities.
  \end{question}
  This problem was posed earlier in \cite{Schippers_confinv_Jyv}, before
  the discovery of Theorem \ref{th:generalized_Nehari}.
  The answer to this problem will certainly involve the generalized
  Nehari invariants (\ref{eq:contour_integral}) and Jenkins' method of quadratic differentials.
\end{subsection}
\begin{subsection}{Schiffer comparison theory of domains}   \label{se:Schiffer_comparison}
 Another form of the comparison moduli space idea is Schiffer's comparison theory of
 nested domains, where one compares the kernel functions of two
 domains $E$ and $\mathcal{E}$ where $E \subset \mathcal{E}$.
 The ideas  go back to his foundational paper
 with S. Bergman \cite{BergmanSchiffer} and his earlier paper \cite{Schiffer_Duke}.  Those papers deal
 with the special case $\mathbb{C}$; in this case the importance of
 the role of the outer
 domain is not obvious.  The comparison theory of kernels on nested domains
 was principally investigated by Schiffer \cite{Schiffer_collected_1,Schiffer_collected_2},
 so we will attach his name to it.  Schiffer's theory reached its final
 published form in his appendix to the book of R. Courant
 \cite{Courant_Schiffer}.  This
 section draws mainly on that source.  We add some simple but important observations on conformal
 invariance.

 Let $E$ be a planar domain with Green's function $g$.  Define
 two kernel functions as follows.
 \begin{definition} The Bergman kernel of a domain $E$ with Green's function $g$ is
 \begin{equation*} \label{eq:K-kernel}
  K(w,z) = - \frac{2}{\pi} \frac{\partial^2 g}{\partial \bar{w}
  \partial z}.
 \end{equation*}
 The Schiffer kernel of $E$ is
 \begin{equation*}  \label{eq:L-kernel}
  L(w,z) = - \frac{2}{\pi}\frac{\partial^2 g}{\partial
  w \partial z}.
 \end{equation*}
 \end{definition}
  The first appearance of the $L$-kernel is to our knowledge \cite{Schiffer_Duke},
  hence we call it the Schiffer kernel.
 Clearly
 \begin{equation*} \label{eq:Lsymmetric}
  L(z,w)=L(w,z)
 \end{equation*}
 and
 \begin{equation*} \label{eq:Ksesquisymmetric}
 K(z,w)= \overline{K(w,z)}.
 \end{equation*}
 These kernels are conformally invariant in the sense that if $g:E \rightarrow E'$
 is a conformal bijection then
 \[  L_{E'}(g(w),g(z)) g'(w) g'(z) = L_E(w,z) \ \ \text{and}
   \ \  K_{E'}(g(w),g(z)) \overline{g'(w)}{g'(z)} = K_E(w,z).   \]
 Strictly speaking it is best to view the kernel functions as differentials $L(w,z)dwdz$
 and $K(w,z) d\bar{w} d{z}$.

 We will define the Bergman and Schiffer kernels of the plane as follows.
 \begin{definition}
  The Bergman kernel of $\mathbb{C}$ is defined to be $K(z,w)=0$.
  The Schiffer kernel of $\mathbb{C}$ is
  \[  L_{\mathbb{C}}(z,w)= \frac{1}{\pi} \frac{1}{(z-w)^2}.  \]
 \end{definition}
 These definitions can be motivated by computing the Schiffer kernel
 of the disk of radius $r$ and letting $r \rightarrow \infty$.  Another motivation
 is that with that definition, the identities and inequalities continue to hold in the case that $\mathcal{E}=\mathbb{C}$.

  Next we give two important properties of the kernel functions.
  Denote the Bergman space of a domain $E$ by
  \[  A_1^2(E) = \left\{ h:E \rightarrow \mathbb{C} \ \text{holomorphic} \ \,:\, \iint_E |h|^2 dA \right\}  < \infty. \]
  \begin{remark}  \label{re:Bergman_is_differentials}  The Bergman space is best viewed as a space of
   one-differentials, but for now we view it as a function space to be consistent with the majority of the function
   theory literature.
  \end{remark}
 \begin{proposition} \label{pr:K_and_Lreproduces}
  Let $f \in A^2_1(E)$.  Then
  \begin{equation*}
   \iint_{E} K(\zeta,\eta)\, f(\eta)\, dA_\eta  =
   f(\zeta)
  \end{equation*} and
  \begin{equation*}
   \iint_{E} L(\zeta,\eta)\, \overline{f(\eta)}\,
   dA_\eta  =  0.
  \end{equation*}
 \end{proposition}

 Let $E$, $E^\mathsf{c}$ and $\mathcal{E}$ be
 domains bounded by finitely many smooth curves, satisfying $E \subset
 \mathcal{E}$ and $E^\mathsf{c} \subset \mathcal{E}$, in such a way that
 $\mathcal{E}$ is the union of $E$ and $E^\mathsf{c}$ together with their
 shared boundary curves.  Note that $E^\mathsf{c}$ is not the literal complement
  in $\mathcal{E}$ but rather its interior.  Denote by $E \bigsqcup E^\mathsf{c}$
the disjoint union of $E$ and $E^\mathsf{c}$. The inner products on each
 domain will be denoted $(,)_{E}$ etc.

  Schiffer considered integral operators naturally associated with
 a configuration of nested domains.  Let $K$, $K^\mathsf{c}$ and $\mathcal{K}$ be
 the Bergman kernels of $E$, $E^\mathsf{c}$ and $\mathcal{E}$ respectively.
 Similarly, let $L$, $L^\mathsf{c}$ and $\mathcal{L}$ be the Schiffer
 L-kernels associated with $E$, $E^\mathsf{c}$ and $\mathcal{E}$.
 \begin{definition}  \label{de:Tone_operator}
 \begin{align*}
  T^1_{E,\mathcal{E}}: A^2_1(E) & \rightarrow A^2_1(\mathcal{E}) \\
  f & \mapsto \iint_E \mathcal{K}(\cdot,w)f(w)\,dA_w
 \end{align*}
 \end{definition}

 Next we consider the operator associated with the Schiffer
 $L$-kernel.
 \begin{definition} \label{de:Ttwo_operator}
  \begin{align*}
   T^2_{E,\mathcal{E}}: \overline{A^2_1(E)} & \rightarrow A^2_1(E \sqcup E^\mathsf{c} ) \\
   \overline{f} & \mapsto \iint_E \mathcal{L}(\cdot,w) \overline{f(w)} \,dA_w.
  \end{align*}
 \end{definition}
 Here the inner product on $A^2_1(E \sqcup E^\mathsf{c} )$ is of course
 \[  (h_1,h_2)_{E \sqcup E^\mathsf{c}}= \iint_{E \sqcup E^\mathsf{c} }
    h_1\overline{h_2}.  \]
 \begin{remark} Schiffer defines the operator $T^2_{E,\mathcal{E}}$ as a map from
  $A^2_1(E)$ into $A^2_1(E \sqcup E^\mathsf{c})$.  This is not convenient, since with that
  convention the operator is complex anti-linear.
 \end{remark}

 There is a
 relation between $T^1_{E,\mathcal{E}}$ and $T^2_{E,\mathcal{E}}$.
 \begin{theorem} \label{th:T2T1relation} For all $f,g \in \overline{A^2_1(E)}$,
 \[  (T^2_{E,\mathcal{E}} g,T^2_{E,\mathcal{E}} f)_{E \sqcup E^\mathsf{c}} = (g,f)_E - (T^1_{E,\mathcal{E}}g,T^1_{E,\mathcal{E}}
   f)_{E \sqcup E^\mathsf{c} }.  \]
 Thus
  \[  {T^2_{E,\mathcal{E}}}^\dagger T^2_{E,\mathcal{E}} = I - {T^1_{E,\mathcal{E}}}^\dagger T^1_{E,\mathcal{E}}.  \]
 \end{theorem}

 Theorem \ref{th:T2T1relation} immediately implies the {\it{bounded Grunsky inequalities}} (in Schiffer's form) for the $T$ operator.
 \begin{corollary} \label{co:bounded_Grunsky}
  For $f \in A^2_1(E)$ one has that
  \[  \|T^2_{E,\mathcal{E}}f \|^2_E = \|f\|^2_E -  \|T^1_{E,\mathcal{E}} f\|^2_{E \cup E^\mathsf{c} } -
    \|T^2_{E,\mathcal{E}} f\|_{E^\mathsf{c} }^2 . \]
  In particular,
  \[  \| T^2_{E,\mathcal{E}} f \|_{E^\mathsf{c}} \leq \|f \|_E.  \]
 \end{corollary}
 \begin{remark} \label{re:Grunsky} In the case that $E$ is simply connected,
  the classical Grunsky inequalities can be obtained by setting $E=g(\disk)$
  for some conformal map $g$, $\mathcal{E}=\mathbb{C}$
  and $f$ to be a
  polynomial in $1/z$ in the above corollary. This can be found in Bergman-Schiffer \cite[Section 9]{BergmanSchiffer}.

 \end{remark}

In Schiffer's comparison theory discussed above, one can identify a comparison moduli space as follows.
 We set $\mathcal{C} = \{ (\mathcal{E},E,h) \,:\, h \in A_1^2(E) \}$, where as before $E \subseteq \mathcal{E}$
 are planar domains.
 Let $g:\mathcal{E} \rightarrow \mathcal{E}'$ be a
 bijective holomorphic map. Define the pull-back by
 \begin{eqnarray*} \label{eq:Bergman_pullback_definition}
  g^* : A^2_1(E') & \rightarrow A^2_1(E) \\
  h & \mapsto h \circ g \cdot g'
 \end{eqnarray*}
 (cf Remark \ref{re:Bergman_is_differentials}) and similarly $g^* \overline{h} = \overline{h \circ g \cdot g'}$.
 Then we have
 \begin{proposition} The operators $T_1$ and $T_2$ are conformally
  invariant in the following sense$:$  If $g$ is a one-to-one
  holomorphic map from $E$ onto $E'$, then for all $h \in A^2_1(E')$
  \[  T^1_{E,\mathcal{E}}\, (g^* h) = g^*T^1_{E',\mathcal{E}'}\, (h)   \ \ \ \ \text{and} \ \ \ \
    T^2_{E,\mathcal{E}}\, (g^* \overline{h}) = g^* T^2_{E',\mathcal{E}'}\, (\overline{h}). \]
 \end{proposition}
 Thus if we define
 \[  I^1(\mathcal{E},E,h) = \| T^1_{E,\mathcal{E}} h \|_{A_1^2(E)}  \ \ \ \text{and}
  \ \ \ I^2(\mathcal{E},E,h) = \| T^2_{E,\mathcal{E}} \overline{h} \|_{A_1^2(E)}  \]
 we obtain conformal invariants on $\mathfrak{C}$.  Similarly the norms of $T^i_{E,\mathcal{E}}$ over $A_1^2(E^\mathsf{c})$ are conformally invariant.

 Remark \ref{re:Grunsky} shows that it is possible to derive coefficient inequalities from Corollary
 \ref{co:bounded_Grunsky}.  We summarize the
 familiar procedure \cite{Durenbook,Pommerenkebook} in this notation.  Setting $\mathcal{E}=\mathbb{C}$ and
 $E=g(\mathcal{\mathbb{D}})$ for a univalent function $g:\disk \rightarrow \mathbb{C}$, it is possible to write the operators
 $T^2_{g(\mathbb{D}),\mathbb{C}}$ and $T^1_{g(\mathbb{D}),\mathbb{C}}$ in terms of the
 coefficients of $g$.  Assuming that $g(0)=0$, for every choice of test function $h$ which is polynomial in $1/z$, one obtains a distinct inequality from Corollary
 \ref{co:bounded_Grunsky}. The same can be done for bounded univalent functions
 by replacing $\mathbb{C}$ with $\mathbb{D}$.
 The inequalities derived in this manner are conformally invariant.

 \begin{remark}  The similarity to the Nehari-type invariants is of course not a
  coincidence.
 \end{remark}

 Some generalizations of Schiffer's comparison theory to Riemann surfaces can be found in
 M. Schiffer and D. Spencer \cite{Schiffer_Spencer}.
\end{subsection}

\begin{subsection}{Fredholm determinant and Fredholm eigenvalues}
The connection between the so-called Fredholm eigenvalues and Grunsky matrices was first observed by Schiffer, see e.g. \cite{Schiffer_Pol}.  One of Schiffer's accomplishments was to place Fredholm's real integral equations of potential
  theory in the complex setting.  We
   draw on his insights freely here without justification, and refer the interested reader to \cite{Schiffer_rend_Fredholm,Schiffer_Pol}.

   Given a  domain $D$ in the plane, bounded by finitely many $C^3$ Jordan curves, consider
   the operator defined for $f\in A_{1}^{2}(D)$ by

\begin{equation}
  T_D f := (T^2_{D,\mathbb{C}} f)|_D,
\end{equation}
where $T^2$ is the operator defined in Definition \ref{de:Ttwo_operator}. The operator $T_D$ is a bounded operator from $\overline{A_{1}^{2}(D)}$ to $A_{1}^{2}(D)$.  We note that this can even be extended to disconnected domains \cite{Schiffer_rend_Fredholm}.

Moreover $T_D$ is a Hilbert-Schmidt operator. This yields that $T_D\overline{T_D}$ is a positive trace class operator, to which the spectral theorem is applicable ($T_D\overline{T_D}$ is compact and self-adjoint). According to this theorem, there is an orthogonal basis $f_n$ of eigenfunctions of $T_D\overline{T_D}$ with non-negative eigenvalues $\lambda_n^2$.
  The non-zero $\lambda_n$ are called the {\it{Fredholm eigenvalues}} of $D$. The reader should however observe that this is not standard in the literature, as some authors refer to $1/\lambda_n$ as the Fredholm eigenvalues. Now since $T_D\overline{T_D}$ is a trace class operator, one can define the  Fredholm determinant  of the domain $D$ as

  \begin{equation}\label{defn:fred determinant}
   \Delta_D:= \det(I- T_D\overline{T_D})= \prod_{n} (1-\lambda_n ^2).
  \end{equation}

  Using variational techniques, Schiffer showed that the Fredholm determinant is conformally invariant in
  the sense of (\ref{eq:conf_inv_nested_surfaces}).  More precisely, let $D_1,\ldots,D_n$ be a collection of bounded simply connected
  domains in $D$, each bounded by a $C^3$ Jordan curve.  Let $E$ be the interior of the complement of $D_1 \cup \cdots \cup D_n$
  in $\mathbb{C}$.  Fix an $l \in \{1,\ldots,n \}$ and let $\mathcal{E} = \mathbb{C} \backslash D_l$.
  Let $\mathfrak{C} = \{ (\mathcal{E},E) \}$ where $E$ and $\mathcal{E}$ are of the above form.
  In that case, the quantity
  \[  I(\mathcal{E},E) = \frac{\Delta_\mathcal{E}}{\Delta_E}   \]
  is a conformal invariant in the sense of (\ref{eq:conf_inv_nested_surfaces}).
  Observe that this striking invariance property must involve pairs of domains, and therefore is a result
  regarding comparison moduli spaces.

  \begin{remark}
  Although we do not doubt the validity of this result, Schiffer's proof is rather heuristic, since it infers invariance under all conformal maps from invariance under a specific  set of variations.
  \end{remark}

   In the case that $D$ is simply connected, $T_D$ is a version of the Grunsky operator, as we have
    observed above.  By work of L. Takhajan and L.-P. Teo \cite{Takhtajan_Teo_Memoirs}, the Fredholm determinant relates to the Weil-Petersson
    metric.  We will re-visit this in Section \ref{se:Kahler} ahead.
\end{subsection}
\end{section}
\begin{section}{Teichm\"uller space as a comparison moduli space}\label{se:Teich_theory}
\begin{subsection}{Spaces of non-overlapping maps and the rigged moduli space}
 Sets of maps with non-overlapping images are often considered in geometric
 function theory.  In the general case of the sphere with $n$ non-overlapping maps,
 these are sometimes called the ``Goluzin-Lebedev class'' \cite{Grinshpan}.  As we
 saw in Example \ref{ex:rigged_moduli_space_analytic}, sets of non-overlapping conformal
 maps appear in conformal field
 theory.   We will be
 concerned with the case that the mappings are quasiconformally extendible
 and the closures do not intersect.

 Let
 \[ A_1^\infty(\disk) = \{ f:\disk \rightarrow \mathbb{C} \ \ \text{holomorphic} \,:\,
   \sup_{z \in \disk} (1-|z|^2) |f(z)| < \infty \}  \]
 and recall that $A_1^2(\disk)$ is the Bergman space on the disk.
 \begin{definition}
  Let $\Oqc$ denote the set of maps $f:\disk \rightarrow \mathbb{C}$
  such that $f$ is an injective holomorphic map, $f(0)=0$ and $f$ is
  quasiconformally extendible to a neighbourhood of $\mathrm{cl}\disk$.
  Let $\Oqc_{\mathrm{WP}} = \{ f \in \Oqc \,:\, f''/f' \in A_1^2(\disk) \}$.
 \end{definition}
 Here ``$\mathrm{cl}$'' denotes closure, and
 ``WP'' stands for Weil-Petersson.  The terminology and motivation for
 this definition, as well as a review of the literature, will be given
 in Section \ref{se:WeilPetersson} ahead.

 \begin{definition}
  Let $\riem$ be a Riemann surface of finite genus with $n$ ordered punctures
  $(p_1,\ldots,p_n)$.
  $\Oqc(\riem)$ is the set of $n$-tuples of one-to-one
  holomorphic mappings $(f_1,\ldots,f_n)$ such that
  \begin{enumerate}
   \item $f_i:\disk \rightarrow \riem$ are one-to-one,
    holomorphic for $i=1,\ldots,n$
   \item $f_i(0)=p_i$
   \item $f_i$ has a quasiconformal extension to a neighbourhood
    of the closure of $\disk$ for all $i$
   \item $f_i(\mathrm{cl}{\disk}) \cap f_j(\mathrm{cl}{\disk})$ is
    empty whenever $i \neq j$.
  \end{enumerate}
 \end{definition}
 This defines a configuration space of mappings in the sense of Section \ref{se:moduli_spaces_mappings}.
 We will call such an $n$-tuple a ``rigging'' of $\riem$.

 Observe also that $(f_1,\ldots,f_n) \in \Oqc(\riem)$ if and only if there is at least one
 $n$-tuple of local coordinates $(\zeta_1,\ldots,\zeta_n)$ on open neighbourhoods
 $U_i$ of $\text{cl} f(\disk)$ such that $\zeta_i \circ f_i \in \Oqc$ for $i=1,\ldots,n$.
 We call such a chart an $n$-chart. If the condition holds for one $n$ chart, then
 it holds for every $n$ chart so that each $U_i$ contains $\text{cl} f_i(\disk) = f_i(\text{cl}{\disk})$.

 Radnell and Schippers showed that the set of riggings is a complex Banach manifold.
 \begin{theorem}[\cite{RSnonoverlapping}]  \label{th:riggings_complex_Banach}
  Let $\riem$ be a Riemann surface with $n$ punctures.  $\Oqc(\riem)$ is a
  complex Banach manifold  locally
 modelled on $\bigoplus^n \left(\mathbb{C} \oplus A_1^\infty(\disk)\right)$.
 \end{theorem}
 The charts
 of this manifold are
 \[  (f_1,\ldots,f_n) \mapsto \left( (\zeta_1 \circ f_1)'(0),\frac{(\zeta_1 \circ f_1)''}{(\zeta_1 \circ f_1)'}, \ldots, (\zeta_n \circ f_n)'(0),\frac{(\zeta_n \circ f_n)''}{(\zeta_n \circ f_n)'} \right)  \]
 for some $n$-chart $(\zeta_1,\ldots,\zeta_n)$.

 Similarly we have the Weil-Petersson class riggings.
 \begin{definition}
  Let $\riem$ be a Riemann surface with $n$ punctures.  $\Oqc_{\mathrm{WP}}(\riem)$ is the
  set of $(f_1,\ldots,f_n) \in \Oqc(\riem)$ such that there is an $n$-chart
  such that $\zeta_i \circ f_i \in \Oqc_{\mathrm{WP}}$.
 \end{definition}
 Again, this defines a configuration space of mappings in the sense of (\ref{eq:map_invariance}).

 Analogously to Theorem \ref{th:riggings_complex_Banach}, Radnell, Schippers and Staubach
 showed that $\Oqc_{\mathrm{WP}}(\riem)$ is a complex Hilbert manifold.
 \begin{theorem}[\cite{RSS_Filbert1}]
  Let $\riem$ be a Riemann surface of finite genus with $n$ ordered punctures.
  $\Oqc_{\mathrm{WP}}(\riem)$ is a complex Hilbert manifold locally modelled on
  $\bigoplus^n \left(\mathbb{C} \oplus A_1^2(\disk)\right)$.
 \end{theorem}
 The charts are obtained in the same way as for $\Oqc(\riem)$.

 We now define the rigged moduli space, which appears in conformal field theory
 and is closely related to the theory of vertex operator algebras \cite{Huang}.
 \begin{definition}
  The{\it{ quasiconformally rigged moduli space}}  is the set of equivalence
  classes
  \[ \widetilde{\mathcal{M}}(g,n) = \{( R, \mathbf{f}) \}/\sim  \]
  where $R$ is a Riemann surface of genus $g$ with $n$ punctures, $\mathbf{f} \in \Oqc(R)$
  and the equivalence relation is defined by $(R,\mathbf{f}) \sim (S,\mathbf{g})$ if and only
  if there is a biholomorphism $\sigma:R \rightarrow S$ such that $\sigma \circ \mathbf{f} =
  \mathbf{g}$.
  The Weil-Petersson class rigged moduli space is the moduli space $ \widetilde{\mathcal{M}}_{\mathrm{WP}}(g,n)$
  defined as above with $\Oqc(R)$ replaced with $\Oqc_{\mathrm{WP}}(R)$.
 \end{definition}
 By $\sigma \circ \mathbf{f}=\mathbf{g}$ we mean that $\sigma \circ f_i = g_i$ for
 $i=1,\ldots,n$ where $\mathbf{f}=(f_1,\ldots,f_n)$ and $\mathbf{g}=(g_1,\ldots,g_n)$.

 In two-dimensional conformal field theory, various
 analytic choices for the riggings exist.
 We will see in the next section that the choices
  above lead to deep connections with Teichm\"uller theory.

 A problem of immediate interest is the following:
 \begin{question}  \label{qu:extend_to_rigged}
  Extend the invariants (\ref{eq:contour_integral}) to functions on $\widetilde{\mathcal{M}}(g,n)$;
  that is, to configuration spaces of maps into Riemann surfaces whose images are bounded by quasicircles.
 \end{question}
 This problem requires some clarification.  The extension to Riemann surfaces was already
 accomplished by Nehari \cite{Nehari_some_inequalities} for harmonic functions, and the methods of Schippers
 \cite{Schippers_Israel_Journal,Schippers_quadratic_expo} for quadratic differentials extend without difficulty to fairly arbitrary Riemann surfaces.
  The main issue is analytic, and amounts to extending the functionals to Riemann surfaces
  bordered by quasicircles.  However, although the
 functionals can be extended continuously from maps with analytic boundaries to roughly
 bounded ones by exhaustion, it is highly desirable to have a natural definition which does not
 require special pleading.
 This could perhaps be done by replacing the contour integrals with appropriate
 integrals with respect to harmonic measure on the boundary of each quasicircle.

 A closely related problem is the following.
 Since the extremals of the functionals map onto $R$ minus trajectories
 of the quadratic differential, we ask
 \begin{question}  \label{qu:rigged_boundary}
  Can the rigged moduli space be endowed with a boundary which includes riggings that are
  maps onto $R$ minus trajectories of the quadratic differential?  This must
  be done in such
  a way that the functionals extend continuously to this boundary.
 \end{question}
\end{subsection}
\begin{subsection}{A brief primer on quasiconformal Teichm\"uller theory}
 \label{se:qc_Teich_primer}  We give an overview of quasiconformal Teichm\"uller theory; for details
 see O. Lehto \cite{Lehto}, S. Nag \cite{Nagbook}, or J. Hubbard \cite{Hubbard_Book}.
  First we define the Teichm\"uller space of a Riemann surface.  Let
 $\riem$ be a Riemann surface whose universal cover is the disk.  Given
 another Riemann surface $\riem_1$, a {\it marking} of this surface by $\riem$
 is a quasiconformal map $f:\riem \rightarrow \riem_1$.  Note that in particular
 this implies that $\riem_1$ is quasiconformally equivalent to $\riem$ (and
 thus in particular they are homeomorphic).  We will denote marked Riemann
 surfaces by triples $(\riem,f,\riem_1)$.

 We say that two quasiconformally equivalent marked Riemann surfaces are Teichm\"uller
 equivalent
 \[  (\riem,f_1,\riem_1) \sim (\riem,f_2,\riem_2)   \]
 if there is a biholomorphism $\sigma:\riem_1 \rightarrow \riem_2$ such that
 $f_2^{-1} \circ \sigma \circ f_1$ is homotopic to the identity rel boundary.
 The term ``rel boundary'' means that the homotopy fixes the boundary pointwise.
 Making the meaning of boundary and homotopy rel boundary precise
 requires some effort in terms of
 the lift to the universal cover; we refer the
 readers to \cite{Hubbard_Book,Lehto} for a thorough treatment.  The
 Teichm\"uller space can now be defined as follows.
 \begin{definition}  \label{de:Teich_space} Let $\riem$ be a fixed Riemann surface.
  The {\it{Teichm\"uller space}} of a Riemann surface $\riem$ is the set of
  equivalence classes
  \[  T(\riem) = \{ (\riem,f,\riem_1) \}/\sim  \]
  where $f:\riem \rightarrow \riem_1$ is a quasiconformal marking of a Riemann
  surface $\riem_1$ and $\sim$ denotes Teichm\"uller equivalence.
  Denote equivalence classes by $[\riem,f,\riem_1]$.
 \end{definition}
 An equivalent definition is given in terms of Beltrami differentials.
 A $(-1,1)$-differential is one which is given in a local coordinate $z$ by
 \[  \mu = h(z) \frac{d\bar{z}}{dz}  \]
 for a Lebesgue-measurable complex-valued function $h$,
 and which transforms under a change of coordinate $z=g(w)$ by
 \[  \tilde{h}(w) = h(g(w)) \frac{\overline{g'(w)}}{g'(w)}; \ \ \ \ \text{that is}  \ \ \   h(z) \frac{d\bar{z}}{dz} = \tilde{h}(w) \frac{d\bar{w}}{dw}.  \]
 It is evident from the transformation rule that the essential supremum of a Beltrami
 differential is well-defined.  For a Riemann surface $\riem$ define
 \begin{equation} \label{eq:Beltrami_definition}
  L^\infty_{-1,1}(\riem)_1 = \{\text{ measurable (-1,1) differentials}\,\, \mu\,\, \text{on}\,\, \riem \,:\,
       \| \mu\|_\infty <1 \}.
 \end{equation}
 Given any $\mu \in L^\infty_{-1,1}(\riem)_1$ there is a quasiconformal solution $f$ to the
 Beltrami differential equation
 \begin{equation} \label{eq:Belt_of_function_def}
   \mu = \frac{\overline{\partial} f}{\partial f}
 \end{equation}
 which is unique up to post-composition by a conformal map.  We call $\mu$ in (\ref{eq:Belt_of_function_def})
 the Beltrami differential of $f$.
 We say two Beltrami differentials $\mu,\nu$ are Teichm\"uller equivalent
 $\mu \sim \nu$ iff the corresponding
 solutions $[\riem,f_\mu,\riem_\mu]$ and $[\riem,f_\nu,\riem_\nu]$ are
 Teichm\"uller equivalent in the sense of Definition \ref{de:Teich_space}.
 The Teichm\"uller space can thus be identified with
 \begin{equation} \label{eq:Beltrami_Teichmuller_defn}
  T(\riem) = \{ \mu \in L^\infty_{-1,1}(\riem)_1 \}/ \sim
 \end{equation}
 via the map
 \[  [\riem,f,\riem_1] \mapsto [\mu].   \]

 In this paper, we restrict our attention to two possible kinds of Riemann
 surfaces: those of genus $g$ with $n$ punctures, and those of genus $g$ and
 $n$ boundary curves homeomorphic to the circle.  In the border case we
 assume that the boundary curves are borders in the sense of L. Ahlfors and L. Sario
 \cite{Ahlfors_Sario}.  We refer to these surfaces as punctured Riemann surfaces
 of type $(g,n)$ or bordered surfaces of type $(g,n)$.

 We now describe the universal Teichm\"uller space, and in doing so fill out the details of Examples \ref{ex:univ_Teich_nested} and \ref{ex:univ_Teich_mapping}.   Let $\disk^* = \{ z\,:\, |z|>1 \} \cup \{\infty \}$.
 The universal Teichm\"uller space $T(\disk^*)$ can be represented as follows.  For
 a given quasiconformal marking $(\disk^*,f,\Omega)$ we let $\hat{f}:\sphere \rightarrow \sphere$ be the quasiconformal map of the Riemann sphere which is conformal on $\disk$ and whose Beltrami
 differential equals that of $f$ on $\disk^*$.  That is,
 \[  \frac{f_{\bar{z}}}{f_z} = \frac{\hat{f}_{\bar{z}}}{\hat{f}_z}  \ \ \text{a.e. on} \ \disk^*.  \]
 The resulting map $\hat{f}$ is unique up to post-composition with a M\"obius transformation
 (and is traditionally specified through three normalizations).  Stated in terms of the
 rigged moduli space, we have that the map
 \begin{align*}
  \Phi:T(\disk^*) & \rightarrow \mathcal{M}(0,1) \\
  [\disk^*,f,\Omega] & \mapsto [\sphere,\hat{f}]
 \end{align*}
 is a bijection.

Given $\mu \in L^\infty_{-1,1}(\disk^*)$ as above,
let $\hat{f}_\mu:\sphere \rightarrow \sphere$ be a quasiconformal map with Beltrami differential equal to $\mu$ in $\disk^*$ and to $0$ in $\disk$.  However, now we uniquely specify $\hat{f}_\mu$ with the normalizations $\hat{f}_\mu(0)=0$, $\hat{f_\mu}'(0)=1$ and $\hat{f}_\mu(\infty)=\infty$.
Denote
\[  f_\mu = \left.  \hat{f}_\mu  \right|_{\disk}  \]
which is a conformal map; this map $f_\mu$ is independent of the choice of representative
in the Teichm\"uller equivalence class.  It can be shown that two Beltrami differentials $\mu$ and $\nu$ on $\disk^*$ represent the same point of $T(\disk^*)$ if and only if $f_\mu=f_\nu$ \cite{Lehto}.  Thus, $T(\disk^*)$ can be identified with the moduli space of Example \ref{ex:univ_Teich_mapping} up to a change of normalization, by applying the uniformization theorem to identify $R$ with $\sphere$.

Define the spaces of abelian differentials
\[  A_1^\infty(\disk) = \{ \alpha(z) dz  \ \text{ holomorphic on } \disk \,:\,  \sup_{z \in \disk}
    (1-|z|^2) |\alpha(z)| <\infty  \} \]
and quadratic differentials
\[  A_2^\infty(\disk) = \{ Q(z)dz^2 \ \text{ holomorphic on } \disk \,:\,
   \sup_{z \in \disk} (1-|z|^2)^2 |Q(z)| < \infty \}.  \]

  Denote the Schwarzian derivative of a conformal map $f$ by
  \[ \mathcal{S}(f) = \frac{f'''}{f'} - \frac{3}{2}
   \left(\frac{f''}{f'} \right)^2.  \]
The Bers embedding is defined by
\begin{align*}  \label{eq:Bers_embedding}
 \beta : T(\disk^*) & \rightarrow A_2^\infty(\disk) \\
 [\mu] & \mapsto \mathcal{S}(f_\mu) dz^2.
\end{align*}
where we use the Beltrami differential model of $T(\disk^*)$.  That is, $\beta = \mathcal{S}
  \circ \Phi$.  It is a classical
theorem of Bers that $\beta(T(\disk^*))$ is an open subset of  $A_2^\infty(\disk)$;
in fact this is a homeomorphism with respect to the so-called Teichm\"uller metric.
Thus $T(\disk^*)$ can be given a complex structure from $A_2^\infty(\disk)$.

 Next, we recall
 another model of the universal Teichm\"uller space $T(\disk^*)$ in terms
 of quasisymmetries.
\begin{definition}\label{quasisymmetriccircle}
An orientation-preserving homeomorphism $h$ of $\mathbb{S}^1$ is called a \emph{quasisymmetric mapping}, iff there is a constant $k>0$, such that for every $\alpha$ and every $\beta$ not equal to a multiple of $2\pi$, the inequality
 \[  \frac{1}{k} \leq \left| \frac{h(e^{i(\alpha+\beta)})-h(e^{i\alpha})}{h(e^{i\alpha})-h(e^{i(\alpha-\beta)})} \right|
    \leq k \]
holds. Let $\qs(\mathbb{S}^1)$ denote the set of quasisymmetric maps from $\mathbb{S}^1$ to $\mathbb{S}^1$.
\end{definition}
The boundary values of a quasiconformal map are in general quasisymmetries \cite{Lehto}.
By a classical result due to Beurling and Ahlfors \cite{Lehto}, any quasisymmetry $h$
has a quasiconformal extension to $\disk^*$.  Another quasiconformal
extension $E(h)$ was given by A. Douady and C. Earle \cite{Nagbook} with the property that
$E(T \circ h \circ S)= T \circ E(h) \circ S$ for any disk automorphisms $T$ and $S$.

  Two Beltrami differentials
 $\mu$ and $\nu$ are Teichm\"uller equivalent if and only if the quasiconformal solutions
 $f^\mu:\disk^* \rightarrow \disk^*$
 and $f^\nu:\disk^* \rightarrow \disk^*$ to the Beltrami equation on $\disk^*$ are equal on $\mathbb{S}^1$
 up to post-composition by a M\"obius transformation of $\mathbb{S}^1$.  That is,
 if and only if there
 is a M\"obius transformation $T:\disk^* \rightarrow \disk^*$ such that $\left.
 T \circ f^\mu \right|_{\mathbb{S}^1} = \left. f^\nu \right|_{\mathbb{S}^1}$.
 Thus we may identify
 \begin{equation} \label{eq:qs_Teich_correspondence}
  T(\disk^*) \cong \qs(\mathbb{S}^1) / \text{M\"ob}(\mathbb{S}^1).
 \end{equation}
 The identification is given by $[\mu] \mapsto \left. f^\mu \right|_{\mathbb{S}^1}$.
 The inverse is obtained by applying the Ahlfors-Beurling extension theorem.

 It is an important result that $\qs(\mathbb{S}^1)$ and $\qs(\mathbb{S}^1)/\text{M\"ob}(\mathbb{S}^1)$
 are groups under composition.  Although the universal Teichm\"uller space has a topological
 structure (determined for example by the Teichm\"uller metric), it is not a topological group,
 since while right composition is continuous, left composition is not.

 A remarkable property of the universal Teichm\"uller space is that it contains
 the Teichm\"uller spaces of all surfaces whose universal cover is the disk.
 This fact is obtained through the representation of surfaces by quotients
 of the disk by Fuchsian groups.  Teichm\"uller theory can also be viewed as the
 deformation theory of Fuchsian groups \cite{Hubbard_Book,Lehto,Nagbook}.  In this paper we take instead the equivalent
 view of Teichm\"uller theory as the space of deformations of complex structures.
\end{subsection}
\begin{subsection}{The Teichm\"uller/rigged moduli space correspondence}
 We saw in Examples \ref{ex:univ_Teich_nested} and \ref{ex:univ_Teich_mapping} that
 the universal Teichm\"uller space is a comparison moduli space.  In fact, by work
 of Radnell and Schippers \cite{RadnellSchippers_monster}, this holds (up to a
 $\mathbb{Z}^n$ action) for Teichm\"uller
 spaces of Riemann surfaces with more boundary curves and higher genus.
 More precisely, the quotient of the
  Teichm\"uller space of genus $g$ surfaces bordered by $n$ closed curves by a
  $\mathbb{Z}^n$ action is the rigged moduli space $\widetilde{\mathcal{M}}(g,n)$.
  We outline these results here.

 The case of $T(\disk^*)$ as the moduli space of maps into $\sphere$ is special,
 since all Riemann surfaces which are homeomorphic to
 the sphere are biholomorphic to the sphere.   We now deal with the
 case that $g$ is non-zero.
 Given a bordered Riemann surface $\riem^B$ of type $(g,n)$ for $g \geq 1$ and $n \geq 1$, represent
 it as a subset of a punctured Riemann surface $\riem^P$ of type $(g,n)$, in such a way that
 the boundary $\partial \riem^B$ consists of $n$ quasicircles each encircling one
 puncture. This can always be done using a sewing procedure \cite{RadnellSchippers_monster}.
  Given
 an element $[\riem^B,f,\riem^B_1] \in T(\riem^B)$, let $\hat{f}:\riem^P \rightarrow \riem_1^P$
 be a quasiconformal map such that $\hat{f}$ is conformal on
 $\riem^P \backslash \text{cl} \riem^B$ and equal to $f$ on $\riem^B$.
 Now fix a collection of conformal maps $\tau_i:\disk \rightarrow \riem^P$ onto each of the
 connected components of the complement $\riem^P \backslash \text{cl} \riem^B$.  We then
 have the map
 \begin{align*}
  \Phi: T(\riem^B) & \rightarrow \widetilde{\mathcal{M}}(g,n) \\
  [\riem^B,f,\riem^B_1] & \mapsto \left[ \riem_1^P,\left( \hat{f} \circ \tau_1,\ldots,
  \hat{f} \circ \tau_n \right) \right].
 \end{align*}

 The pure Teichm\"uller modular group $\text{PMod}(\riem^B)$ is the set of equivalence
 classes of quasiconformal maps $\rho:\riem^B \rightarrow \riem^B$ such that $\rho(\partial_i \riem^B) = \partial_i \riem^B$ as a set.  Two such maps $\rho_1,\rho_2$ are equivalent $\rho_1 \sim \rho_2$
 if and only if they are homotopic rel boundary.  The modular group acts
 discontinuously via $[\rho][\riem^B,f,\riem^B_1] = [\riem^B, f \circ \rho^{-1},
 \riem^B_1]$, and for any $[\rho]$ the induced map of $T(\riem^B)$ is a
 biholomorphism.  Let $\text{DB}$ denote the subset of $\text{PMod}(\riem^B)$ generated
 by quasiconformal maps which are the identity on every boundary curves and are homotopic
 to the identity map.  This group can be pictured as the set of elements twisting
 each boundary curve an integer number of times, and is isomorphic to $\mathbb{Z}^n$
 except when $g=0$ and $n=1$ or $n=2$, in which case it is trivial or isomorphic to $\mathbb{Z}$
 respectively.

 We then have the following.
 \begin{theorem}[\cite{RadnellSchippers_monster}]  \label{th:RS_corres_thm} Let $\riem^B$ be a bordered Riemann surface.     The map
  $\Phi: T(\riem^B) \rightarrow \widetilde{\mathcal{M}}(g,n)$
  is a bijection up to a discrete group action.  That is, if $\Phi([\riem^B,f_1,\riem^B_1]=\Phi([\riem^B,f_2,\riem^B_2])$
  then there is a $[\rho] \in \text{DB}$ such that $[\rho][\riem^B,f_1,\riem^B_1] =
  [\riem^B,f_2,\riem^B_2]$.
 \end{theorem}

 Here is a specific example.
 \begin{example}[\cite{RadnellSchippers_annulus}]
  Let $A$ be an annulus of finite modulus.  The Teichm\"uller space $T(A)/\mathbb{Z}$ can
  be identified with the rigged moduli space $\widetilde{\mathcal{M}}(0,2)$.

  Interestingly this has a semigroup structure (it is known as the {\it{Neretin-Segal semigroup}})
  and every element can be decomposed as the product of a quasisymmetry of $\mathbb{S}^1$ and a bounded
  univalent function.
 \end{example}

 This example suggests the following problems, closely related to Problem \ref{qu:extend_to_rigged}:
 \begin{question}  \label{qu:extend_to_TA} Lift the conformally invariant functionals (\ref{eq:contour_integral}) to $T(A)$.
 \end{question}
 Of course this problem extends to Teichm\"uller space.
 \begin{question}  \label{qu:extend_to_Teich} Extend the functionals (\ref{eq:contour_integral}) to
  Teichm\"uller spaces of bordered Riemann surfaces of type $(g,n)$.
 \end{question}
 And again we can ask, similarly to Problem \ref{qu:rigged_boundary}:
 \begin{question}  Does the Teichm\"uller space of bordered surfaces of type $(g,n)$
  have a boundary, which corresponds under $\Phi$ to those riggings which map onto analytic arcs of
  quadratic differentials?
 \end{question}

 \begin{question}  Can Theorem \ref{th:RS_corres_thm} be extended to Teichm\"uller spaces
 of more general surfaces,
  say those with infinite genus and/or infinitely many riggings?  What conditions
  are necessary to make this work?
 \end{question}
 It is likely that a bound on the lengths of closed geodesics is a necessary condition;
 see Section \ref{se:WP_higher_genus} and Problem \ref{qu:extend_to_infinite_genus_WP} ahead.
\end{subsection}
\begin{subsection}{Some applications of the Teichm\"uller space/rigged moduli space correspondence}
 This identification of the moduli spaces appearing in conformal field theory (CFT) and in Teichm\"uller
 theory has far-reaching consequences.  In particular, one can endow the
 rigged moduli space with a complex structure and show that
 the operation of sewing bordered Riemann surfaces together with quasisymmetries
 is holomorphic \cite{RadnellSchippers_monster}.  These are required for the construction of two-dimensional conformal field theory
 from vertex operator algebras.
  A comprehensive review can be found in  \cite{RadnellSchippersStaubach_CFT_review}.
  The sewing technique is also of independent interest in Teichm\"uller theory.

 On the other hand, the correspondence can be applied to transfer structures
 in CFT to Teichm\"uller theory.  For example, one can obtain a fibre structure on Teichm\"uller space as follows.
 Let $\riem^B$ be a fixed Riemann surface.  Define the map
 \begin{align*}
  \mathcal{C}: T(\riem^B) & \rightarrow T(\riem^P) \\
  [\riem^B,f,\riem_1^B] & \mapsto [\riem^P,\hat{f},\riem_1^P]
 \end{align*}
 where $\riem^P$, $\hat{f}$ and $\riem_1^P$ are determined from $[\riem^B,f,\riem^B_1]$
 as above.  Radnell and Schippers showed the following.
 \begin{theorem}[\cite{RadnellSchippers_fiber}] \label{th:fiber_structure}
  Let $\riem^B$ be a bordered surface of type $(g,n)$, and $\riem^P$
 be a punctured surface of type $(g,n)$.   Assume that $\riem^B \subset \riem^P$
 and $\partial \riem^B$ consists of $n$ quasicircles each enclosing a distinct puncture
 in $\riem^P$.
  \begin{enumerate}
   \item $\mathcal{C}$ is a holomorphic map with local holomorphic sections.
    The fibres $\mathcal{C}^{-1}(p)$ for $p \in T(\riem^P)$ are complex
    submanifolds of $T(\riem^B)$.
   \item Let $p \in T(\riem^P)$.  If $(\riem^P,g,\riem_1^P)$ is a representative
   of $p$, then the quotient $\mathcal{C}^{-1}(p) / \text{DB}$ is canonically bijective
   with $\Oqc(\riem_1^P)$.
   \item The bijection in (2) is a biholomorphism.
  \end{enumerate}
 \end{theorem}
 In particular, this can be used to give holomorphic coordinates on $T(\riem^B)$.
 The proof uses a variational technique of F. Gardiner  \cite{Gardiner}, which was based on an idea
  of Schiffer.  The use of Gardiner-Schiffer variation to address holomorphicity
 issues in the rigged moduli space originates with Radnell \cite{Radnell_thesis}.

 \begin{question}  Are there generalizations of Theorem \ref{th:fiber_structure} for
  Teichm\"uller spaces to more general Riemann surfaces, such
  as those with infinite genus and/or infinite number of boundary curves?
 \end{question}
\end{subsection}
\end{section}
\begin{section}{Weil-Petersson class Teichm\"uller space}   \label{se:WeilPetersson}
\begin{subsection}{The Weil-Petersson class Teichm\"uller space}
 There exist many refinements of Teichm\"uller space, including for example
 the asymptotically conformal Teichm\"uller space \cite{EarleGardinerLakic}, BMO-Teichm\"uller space \cite{CuiZinsmeister},
 and $L^p$
 Teichm\"uller spaces (references below).  The $L^2$ case is usually referred to as the Weil-Petersson class Teichm\"uller
 space, which we motivate in this section.  From now on we abbreviate this as WP-class.

 The Weil-Petersson metric is a metric in the sense of Riemannian/Hermitian geometry, that is,
 an inner product at every tangent space.
 It is based on a representation of the
 tangent space to Teichm\"uller space by a set of quadratic differentials, which we must
 now describe. Fix a
 Riemann surface $\riem$.  Given
 a holomorphic curve $t \mapsto \mu_t$ in the set of Beltrami differentials $L^\infty_{-1,1}(\riem)_1$, the derivative $\frac{d}{dt} \mu$
 is in $L^\infty_{-1,1}(\riem)$.  However, there is an enormous amount of redundancy in this model,
 since the difference between any pair of these might be tangent to the Teichm\"uller equivalence relation, and hence represent the same direction in Teichm\"uller space.

 The following decomposition remedies the problem.  Let $\lambda^2$ denote the
 hyperbolic metric on $\riem$.  Define the harmonic Beltrami differentials by
 \begin{equation} \label{eq:harmonic_Belt_full_def}
  \Omega_{-1,1}(\riem)  = \{ \mu \in L^\infty_{-1,1}(\riem) \,:\, \mu = \alpha/\lambda^2,
 \ \ \ \alpha \text{ a quadratic differential} \}.
 \end{equation}
 It can easily be verified that a quotient of a quadratic differential by the hyperbolic
 metric is a Beltrami differential, by writing them in local coordinates and verifying
 that the quotient satisfies the correct transformation rule.
 It is a classical result that we have the
 Banach space decomposition
 \[  L^\infty_{-1,1}(\riem)_1 = \Omega_{-1,1}(\riem) \oplus \mathcal{N}  \]
 where $\mathcal{N}$ is the set of Beltrami differentials which are tangent to the
 Teichm\"uller equivalence relation. These are the so-called ``infinitesimally trivial
 Beltrami differentials'', which can be characterized precisely \cite{Lehto,Nagbook}.

 Thus, the tangent space to $T(\riem)$ at $[\riem, \text{Id},\riem]$ can be identified
 with $\Omega_{-1,1}(\riem)$.   Furthermore, there is an open subset $U \subset \Omega_{-1,1}(\riem)$ containing $0$ such that
 \begin{align} \label{eq:local_coord_near_id}
  \Psi :U & \rightarrow T(\riem) \nonumber \\
  \mu & \mapsto [\riem,f,\riem_1], \ \ \ \ \ \text{where} \ \ \overline{\partial} f/\partial f = \mu
 \end{align}
 is a biholomorphism from an open neighbourhood of $0$ to an open
 neighbourhood of $0$.  Applying right compositions by quasiconformal
 maps, one obtains a system
 of coordinates on $T(\riem)$ which is compatible with that obtained
 from the Bers embedding $\beta$. This is closely related to the so-called Ahlfors-Weill reflection.
  The tangent spaces at other points can be obtained
 by right composing by quasiconformal maps.  These constructions can also all be lifted to the universal cover and expressed in terms of differentials invariant under Fuchsian groups.

 For compact surfaces or compact surfaces with punctures, the Weil-Petersson metric is defined
 on the tangent space at the identity by the $L^2$ pairing of Beltrami differentials
 \begin{equation} \label{eq:WP_pairing}
  \left< \mu, \nu \right> = \iint_{\riem} \overline{\mu} \nu \,dA_{\mathrm{hyp}}, \ \
   \mu,\nu \in \Omega_{-1,1}(\riem)
 \end{equation}
 where $dA_{\mathrm{hyp}}$ is the hyperbolic area measure \cite{Ahlfors_Kahler}.   The integral (\ref{eq:WP_pairing})
 is finite, because on
 compact surfaces $L^\infty$ Beltrami differentials are also $L^2$. This Hermitian
 inner product can be transferred to tangent spaces at other points
 in the Teichm\"uller space by right composition.  In the compact/compact with punctures case, the Weil-Petersson metric
 has been much studied \cite{Wolpert}. Ahlfors showed that the WP-metric is K\"ahler \cite{Ahlfors_Kahler}
 and later computed its Ricci and scalar curvatures of holomorphic sections, and showed they are negative \cite{Ahlfors_curvature}.

 The Weil-Petersson inner product
 had not been defined on any other Teichm\"uller space until the turn of the millenium.
 The obstacle was that if the Riemann surface is not compact or compact
 with punctures,
 then the $L^2$ pairing need not be finite.  Thus in order to define a Weil-Petersson
 metric one must restrict to a smaller Teichm\"uller space, so that tangent directions
 generate only $L^2$ Beltrami differentials.  Nag and Verjovsky \cite{NagVerjovsky}
 showed that if one restricts to the subset of the universal Teichm\"uller space
 corresponding to diffeomorphisms of $\mathbb{S}^1$, then the WP metric converges on tangent directions.   Equivalently,
 the corresponding representative quasisymmetries are smooth.  However,
 this is only a heuristic principle, and the smoothness assumption is rather
 artificial from the point of view of Teichm\"uller theory.  The correct analytic
 condition giving the largest
 universal Teichm\"uller space on which the WP metric converges is the subject of the next section.
\end{subsection}
\begin{subsection}{Weil-Petersson class universal Teichm\"uller space}
   The investigation of $L^2$ Teichm\"uller spaces
 is due independently to G. Cui \cite{Cui} and L. Takhtajan
  and L.-P. Teo \cite{Takhtajan_Teo_Memoirs}.  H. Guo \cite{GuoHui} and S.-A. Tang \cite{Tang}
  extended some of the results to $L^p$ Teichm\"uller spaces.

The definitions of the WP-class Teichm\"uller space given by  Takhtajan/Teo and Cui are rather different, but equivalent.  Takhtajan/Teo give the definitions
and complex structure in terms of $L^2$ harmonic Beltrami differentials, while
Cui's definition is in terms of an $L^2$ condition on the quadratic differentials
in the image of the Bers embedding.
\begin{remark} \label{re:TT_foliation}
 In fact Takhtajan and Teo's approach defines
 a foliation of $T(\disk^*)$ by right translates of the WP-class Teichm\"uller space.
\end{remark}

Recall the space of differentials
\[  A_1^2(\disk) =   \left\{ \alpha(z)dz  \ \text{ holomorphic on } \disk \,:\,
   \iint_{\disk} |\alpha(z)|^2 < \infty \right\}  \]
and quadratic differentials
\[  A_2^2(\disk) =  \left\{ Q(z)dz^2  \ \text{ holomorphic on } \disk \,:\,
   \iint_{\disk} (1-|z|^2)^2|Q(z)|^2 < \infty \right\}.  \]
The coefficients of the differentials in $A_1^2(\disk)$ are in the Bergman space, and thus
we use the same notation as in Section \ref{se:Schiffer_comparison}  (see Remark \ref{re:Bergman_is_differentials}).

We begin with Cui's definition.
\begin{definition}  \label{de:WP_univ}
The{\it{ $\mathrm{WP}$-class universal Teichm\"uller space}} is defined to be
\[  T_{\mathrm{WP}}(\disk^*) = \{ [\mu] \in T(\disk^*) \,:\, \beta([\mu]) \in  A_2^2(\disk) \}. \]
\end{definition}
Guo, Cui and Takhtajan/Teo also showed that $[\mu] \in T_{\mathrm{WP}}(\disk^*)$ if and
only if $f_\mu''/f_\mu' \in A_1^2(\disk)$.  Thus we have that
\[  [\mu] \in T_{\mathrm{WP}}(\disk^*) \Leftrightarrow f_\mu \in \Oqc_{\mathrm{WP}}(\disk).     \]

Although Takhtajan and Teo used a different definition, they also showed
the equivalence of their definition with Definition \ref{de:WP_univ}, so we attribute the theorems below
jointly.

Takhtajan/Teo and Cui independently showed that
\begin{theorem}[\cite{Cui,Takhtajan_Teo_Memoirs}]  \label{th:disk_inclusion_holomorphic}
$A_2^2(\disk) \subset A_2^\infty(\disk)$
and the inclusion is holomorphic.
\end{theorem}
Combined with the following, we get an analogue
of the Bers embedding for $T_{\mathrm{WP}}(\disk^*)$.
\begin{theorem}[\cite{Cui,Takhtajan_Teo_Memoirs}]  \label{th:Bers_univ_WP}
 $\beta(T_{\mathrm{WP}}(\disk^*))= \beta(T(\disk^*)) \cap A_2^2(\disk)$.  In particular,
 $\beta(T_{\mathrm{WP}}(\disk^*))$ is open.
\end{theorem}
Thus we have that $T_{\mathrm{WP}}(\disk^*)$ has a complex structure inherited from $A_2^2(\disk)$.

Guo and Cui showed that the fact that $\beta([\mu]) \in A_2^2(\disk)$ is equivalent to the existence of a representative
Beltrami differential which is $L^2$ with respect to the hyperbolic metric.  Define for any Riemann surface $\riem$ with
a hyperbolic metric $\lambda^2$
\[  L^2_{\mathrm{hyp}}(\riem) = \left\{ \mu \in L^\infty_{-1,1}(\riem) \,:\, \iint_{\riem} |\mu|^2 dA_{\mathrm{hyp}} < \infty \right\}  \]
where $dA_{\mathrm{hyp}}$ denotes hyperbolic area measure.  In particular
\[  L^2_{\mathrm{hyp}}(\disk^*) = \left\{ \mu \in L^\infty_{-1,1}(\disk^*) \,:\, \iint_{\disk^*} \frac{|\mu|^2}{(1-|z|^2)^2} dA < \infty \right\}.   \]
We then have the following.
\begin{theorem}[\cite{Cui,Takhtajan_Teo_Memoirs}]  \label{th:representative_in_L2}
 $[\mu] \in T_{\mathrm{WP}}(\disk^*)$ if and only if $[\mu]$ has a representative $\mu \in L^2_{\mathrm{hyp}}(\disk^*)$.
\end{theorem}
In fact, Guo and Cui gave stronger results in two directions.
Guo showed such a result for $L^p$ differentials with respect to the hyperbolic metric, $p \geq 1$.
Cui showed that the Douady-Earle extension of the boundary values
of any representative conformal map $f_\mu$ has this property,
and in fact satisfies a stronger integral estimate.    Tang showed that the Douady-Earle
extension is in $L^p$ and that the Bers embedding is holomorphic with respect to the
intersection norm on $L^p \cap L^\infty$.

If $f_\mu \in \Oqc$, then the image of $\mathbb{S}^1$ under the unique homeomorphic
extension of $f_\mu$ to $\text{cl}\disk$ is a quasicircle.  Although there are an astonishing number
of non-trivially equivalent characterizations of quasicircles \cite{Gehring_Hag}, there is no known geometric
characterization of Weil-Petersson class quasicircles.
\begin{question}  Characterize quasicircles of the form $f_\mu(\mathbb{S}^1)$ for $[\mu] \in T_{\mathrm{WP}}(\disk)$ (equivalently, for $f_\mu \in \Oqc_{\mathrm{\mathrm{WP}}}$) via analytic or geometric conditions
on the set itself.
\end{question}
Radnell, Schippers, and Staubach \cite{RSS_WPjump} showed that a Weil-Petersson class
quasicircle is a rectifiable chord-arc curve, but this is unlikely to be sufficient.

We also define the WP-class quasisymmetries as follows.
\begin{definition}
 We say that $\phi \in \qs(\mathbb{S}^1)$ is {\it{Weil-Petersson class}} if its corresponding
 Teichm\"uller space representative is in $T_{\mathrm{WP}}(\disk^*)$.  Denote the set
 of WP-class quasisymmetries by $\qs_{\mathrm{WP}}(\mathbb{S}^1)$.
\end{definition}

Cui and Takhtajan/Teo showed that, like quasisymmetries, these are closed under composition and inverse.
\begin{theorem}[\cite{Cui,Takhtajan_Teo_Memoirs}] $\qs_{\mathrm{WP}}(\mathbb{S}^1)$ is a group.
\end{theorem}
A stronger result was obtained by Takhtajan and Teo:
\begin{theorem}[\cite{Takhtajan_Teo_Memoirs}] $\qs_{\mathrm{WP}}(\mathbb{S}^1)/\text{M\"ob}(\mathbb{S}^1)$
 is a topological group.
\end{theorem}
This is in contrast to $\qs(\mathbb{S}^1)/\text{M\"ob}(\mathbb{S}^1)$, which is not a topological group.
\begin{remark}
 Takhtajan and Teo also showed that $\qs_{\mathrm{WP}}(\mathbb{S}^1)/\mathbb{S}^1$ is a topological group, and it can be identified
 naturally with the WP-class universal Teichm\"uller curve.
\end{remark}
The WP-class quasisymmetries were characterized by Y. Shen as follows, answering a problem
posed by Takhtajan and Teo.
\begin{theorem}[\cite{Shen_characterization}]  \label{th:Shen_characterization}
 Let $\phi:\mathbb{S}^1 \rightarrow \mathbb{S}^1$ be a homeomorphism.  $\phi \in \qs_{\mathrm{WP}}(\mathbb{S}^1)$ if and only if $\phi$ is absolutely continuous and $\log{\phi'}$
 is in the Sobolev space $H^{1/2}(\mathbb{S}^1)$.
\end{theorem}
Further characterizations (e.g. in terms of the composition operator associated with $\phi$) were given by Y. Hu and Y. Shen \cite{HS}.

It was also shown that right composition is biholomorphic.
\begin{theorem}[\cite{Cui,Takhtajan_Teo_Memoirs}]
 Right composition $($mod $\text{M\"ob}(\mathbb{S}^1)$$)$ in $T(\disk^*)$ by a fixed element $h \in \qs_{\mathrm{WP}}(\disk^*)$ is a biholomorphism.
\end{theorem}
This theorem combined with Theorem \ref{th:Bers_univ_WP} was used by Cui to define the Weil-Petersson pairing
on any tangent space, by using the pairing in $A_2^2(\disk)$ at the identity, and then applying the
above theorem to define a right-invariant metric.

Takhtajan and Teo's approach, on the other hand, defined the complex structure
in terms of local charts into the space of harmonic Beltrami differentials, which
are defined by
\[  H_{-1,1}(\disk^*) = \{ \mu = (1-|z|^2)^2 \overline{Q(z)} \,:\,  \mu \in L^2_{\mathrm{hyp}}(\disk^*)  \}.  \]
This is non-trivial as it must be shown that the transition functions of the charts
are biholomorphisms.  As in the classical case, their proof relies on the use of the Ahlfors-Weill reflection.

Their construction also gives a description of the tangent space in terms
of harmonic Beltrami differentials.  By Theorem \ref{th:disk_inclusion_holomorphic} we immediately obtain
\begin{theorem}[\cite{Takhtajan_Teo_Memoirs}] \label{th:TT_tangent_Beltrami}
 \[  L^\infty_{-1,1}(\disk^*) \cap L^2_{\mathrm{hyp}}(\disk^*)
  = H_{-1,1}(\disk^*) \oplus (\mathcal{N} \cap L^2_{\mathrm{hyp}}(\disk^*)).  \]
\end{theorem}
This also shows that the tangent space at $[0]$ to $T_{\mathrm{WP}}(\disk^*)$ can
be identified with $H_{-1,1}(\disk^*)$.
This has finite Weil-Petersson pairing (\ref{eq:WP_pairing}) by definition.
By applying the holomorphic right translation, Takhtajan and Teo
obtain a right invariant inner product at all points in $T_{\mathrm{WP}}(\disk^*)$.

In fact, they extended this complex structure and Hermitian metric to all of $T(\disk^*)$ in the following way.
A neighbourhood of $[0]$ in $T(\disk^*)$ can be obtained under the map $\Psi$ (\ref{eq:local_coord_near_id}),
with $\Omega_{-1,1}(\disk^*)$ replaced by $H_{-1,1}(\disk^*)$.  Using
right composition, the charts patch
together to give a complex structure compatible with the Bers embedding.  However,
the topology and complex structure are not equivalent to the standard one on $T(\disk^*)$.  Indeed,
$T(\disk^*)$ consists of uncountably many disjoint
translates of $T_{\mathrm{WP}}(\disk^*)$, each of which is a connected component of $T(\disk^*)$
  with this new topology.   In their definition,
$T_{\mathrm{WP}}(\disk^*)$ appears as the connected component of the identity.  This accounts for the
convergence of the WP pairing, since at any given point there are far fewer tangent
vectors than there are when $T(\disk^*)$ is given the standard complex structure.  That is, the directions on which the WP pairing diverges are excluded.

It must be emphasized that these constructions require a great deal of
analysis; one cannot simply make small adjustments to classical theorems of Teichm\"uller theory.

Finally, Takhtajan and Teo showed (improving on Nag and Verjovsky's result, which held only in the diffeomorphism
case)
\begin{theorem}[\cite{NagVerjovsky,Takhtajan_Teo_Memoirs}]
 The Weil-Petersson metric on $T(\disk^*)$ $($with the new complex structure and topology$)$ is K\"ahler.
 In particular the Weil-Petersson metric on $T_{\mathrm{WP}}(\disk^*)$ is K\"ahler.
\end{theorem}
In fact they gave different explicit formulas for the K\"ahler potential, and computed the Ricci
curvatures.  We will return to this in Section \ref{se:Kahler}.
\end{subsection}
\begin{subsection}{Higher genus Weil-Petersson class Teichm\"uller spaces}  \label{se:WP_higher_genus}
 It is possible to extend the Weil-Petersson metric to a much wider class
 of surfaces, again by obtaining an $L^2$ theory.
 Radnell, Schippers and Staubach \cite{RSS_Filbert2, RSS_Filbert1,RSS_WP1,RSS_WP2}
 did this for bordered surfaces of type $(g,n)$.  M. Yanagishita \cite{Yanagishita}
 extended the $L^p$ theory of Guo \cite{GuoHui} and Tang \cite{Tang} to
 surfaces satisfying ``Lehner's condition'', which includes bordered
 surfaces of type $(g,n)$. In the $L^2$-setting the two Teichm\"uller spaces are the same as sets, but the constructions of the
 complex structures are rather different.  Yanagishita constructs
 the complex structure from $L^p$ quadratic differentials  under the image of the Bers embedding, following the approach of Cui, Guo
 and Tang.   Radnell, Schippers, and Staubach constructed the complex structure in two equivalent ways:
 using harmonic Beltrami differentials (along the lines of Takhtajan/Teo), and by refining the fiber structure of Theorem \ref{th:fiber_structure}.  As in the classical $L^\infty$ case, the complex structures arising from
 the Bers embedding into quadratic differentials and from harmonic Beltrami differentials should be equivalent,
 but this has not yet been established.\footnotemark

 \footnotetext{See Footnote 3}

 The Weil-Petersson class Teichm\"uller space is specified by boundary behaviour.
 To describe this we need a local coordinate near the boundary.
 \begin{definition}  Let $\riem$ be a bordered surface of type $(g,n)$.
  A collar neighbourhood of a boundary curve $\partial_i \riem$ is a
   doubly connected open set in $\riem$, one of whose boundaries is $\partial_i \riem$
   and the other is an analytic curve in the interior of $\riem$.
   A collar
   chart of a bordered surface $\riem$ is a conformal map $\zeta:U \rightarrow \{1<|z|<r \}$
   for some $r >1$ which extends continuously to $\partial_i \riem$.
 \end{definition}
 It is possible to show that the chart extends to a conformal map of an open neighbourhood
 of $\mathbb{S}^1$ into the double of $\riem$.

 We define the WP-class Teichm\"uller space in two steps.
 \begin{definition} Let $\riem$ and $\riem_1$ be bordered surfaces of type $(g,n)$.  A quasiconformal
  map $f:\riem \to \riem_1$ is called refined if for each pair of boundary curves $\partial_i \riem$, $\partial_j \riem_1$
  such that $f(\partial_i \riem) = \partial_j \riem_1$,
  there are collar charts $\zeta_i, \eta_j$ of $\partial_i \riem$ and $\partial_j \riem_1$ respectively such that $\left. \eta_j \circ f \circ \zeta_i^{-1}
  \right|_{\mathbb{S}^1} \in \qs_{\mathrm{WP}}(\mathbb{S}^1)$.  Denote the set of such quasiconformal
  maps by $\qc_r(\riem)$.
 \end{definition}
 If the condition holds for one collar chart at a boundary $\partial_i \riem$, then it holds for all of them.

 \begin{definition}  \label{de:WP_Teich_general}
  Let $\riem$ be a bordered surface of type $(g,n)$.  The
 WP-{\it {class Teichm\"uller space}} of $\riem$ is
  \[  T(\riem)  = \{ (\riem,f,\riem_1) \,:\, f \in \qc_r(\riem) \} /\sim   \]
  where $\sim$ is the usual Teichm\"uller equivalence.
 \end{definition}
 A different definition was given by Yanagishita \cite{Yanagishita}, for $L^p$ spaces
 for $p \geq 1$.
 It was phrased in terms of Fuchsian groups satisfying a condition he terms
 ``Lehner's condition''.  We will restate Yanagishita's approach in its equivalent
 form on the Riemann surface, for consistency of presentation.  Let $\riem^*$ denote
 the double of the Riemann surface $\riem$ (if $\Gamma$ is the Fuchsian group such that
 $\riem = \disk^*/\Gamma$, then $\riem^* = \disk/\Gamma$).
 \begin{definition}[Lehner's condition]
  A Riemann surface $\riem$ covered by the disk $\disk^*$ satisfies {\it{Lehner's condition}} if the infimum of the
  hyperbolic lengths of the simple closed geodesics is strictly greater than $0$.
 \end{definition}
 \begin{definition}  \label{de:Yanagishita}
  Let $\riem$ be a Riemann surface covered by the disk satisfying Lehner's condition.  The {\it{$p$-integrable
  Teichm\"uller space}} $T^p(\riem)$ is the subset of $T(\riem)$ consisting of elements $[\riem,f,\riem_1]$
  such that there is a representative $(\riem,f,\riem_1)$ such that the Beltrami differential
  of $f$ is in $L^p$ with respect to the hyperbolic metric.
 \end{definition}
 Yanagishita also showed that the $L^p$ representative is given by the Douady-Earle extension
 of the boundary values of the lift to $\disk^*$.
   We will see shortly that it agrees with the
  definition above in the special case of $p=2$ and bordered surfaces of type $(g,n)$.

 Radnell, Schippers and Staubach obtained the following analogue of Theorem \ref{th:representative_in_L2},
 using sewing techniques and the lambda lemma.
 \begin{theorem}[\cite{RSS_WP1}]  \label{th:L2_Beltrami_rep_general}
  Let $\riem$ be a bordered
  Riemann surface of type $(g,n)$.  Then $f \in \qc_r(\riem)$ if and only if
  it is homotopic rel boundary to a quasiconformal map whose Beltrami differential
  is in $L^2_{\mathrm{hyp}}(\riem)$.
 \end{theorem}
 Since bordered surfaces of type $(g,n)$ satisfy
 Lehner's condition, by Theorem \ref{th:L2_Beltrami_rep_general}, the
 Definitions \ref{de:WP_Teich_general} and \ref{de:Yanagishita}  are
 equivalent for $p=2$ and bordered surfaces of this type.

 The two approaches to the complex structure are rather different.  Yanagishita's
 approach involves the following theorem. Let $\riem^*$ denote the double of $\riem$.  Define
 $L^p_{\mathrm{hyp}}(\riem)$
 to be the set of Beltrami differentials on $\riem$ which are
 $L^p$ with respect to the
 hyperbolic area measure.  Following \cite{Yanagishita} denote
 \[Ael^p(\riem) =  L^\infty_{-1,1}(\riem)_1
  \cap L^p_{\mathrm{hyp}}(\riem).   \] The intersection norm $\| \cdot \|_p + \| \cdot \|_\infty$ induces a topology on $T^p(\riem)$.   Furthermore let
  \[  A_2^p(\riem^*) = \left\{  \alpha \ \text{a quadratic differential on} \ \riem^*\,:\,
     \iint_{\riem^*} \lambda^{-2} |\alpha|^p <\infty \right\}.  \]
  (recall that $\lambda^2$ is the hyperbolic metric).
 \begin{theorem}[\cite{Yanagishita}]  For $p \geq 2$, the restriction of the Bers embedding $\beta$ to
  $T^p(\riem)$ is a homeomorphism onto its image in $A_2^p(\riem^*)$ with respect to the
  $Ael^p$ norm.
 \end{theorem}
 This induces a complex structure on $T^p(\riem)$, and in particular on
 $T_{\mathrm{WP}}(\riem)$ for bordered Riemann surfaces of type $(g,n)$.   Furthermore,
 Yanagishita showed that right composition is a biholomorphism with respect to this
 structure, for $p \geq 2$.  Thus, although the tangent space structure is not treated in the
 paper \cite{Yanagishita},
 it is possible to define a WP-pairing on the tangent space at $[0]$ using the $\Gamma$-invariant subspace of $A_2^2(\disk)$ as a model, and then using right composition
 to obtain a right-invariant metric at every point.\\

 Radnell, Schippers, and Staubach gave two other complex structures for
 bordered surfaces of type $(g,n)$, which are equivalent
 to each other.  In \cite{RSS_Filbert2, RSS_Filbert1}, it was shown that the
 fiber structure of Theorem \ref{th:fiber_structure} passes down to $T_{\mathrm{WP}}(\riem)$
 for bordered Riemann surfaces of type $(g,n)$.  That is, one may view $T_{\mathrm{WP}}(\riem)$
 as fibred over $T(\riem^P)$ for a compact surface with punctures $\riem^P$, such that the
 fibres $\mathcal{C}^{-1}(p)$ modulo a discrete group action
 are biholomorphic to $\Oqc_{\mathrm{WP}}(\riem^P_1)$.  This can be used to construct a Hausdorff,
 second countable topology on $T_{\mathrm{WP}}(\riem)$, and a complex Hilbert
 manifold structure \cite{RSS_Filbert2}.
 The advantage of this approach is that it is very flexible and
 constructive, and explicit coordinates can be given in terms of Gardiner-Schiffer
 variation.

 In \cite{RSS_WP2}, Radnell, Schippers, and Staubach showed that this fiber structure is compatible with that obtained
 by harmonic Beltrami differentials and right translation, analogous to both the classical
 case and Takhtajan and Teo's approach on $T_{\mathrm{WP}}(\disk^*)$.  We briefly describe this below,
 as well as the description of the tangent spaces.

 In all of the following
 theorems, $\riem$ is a bordered surface of type $(g,n)$.
 \begin{theorem}[\cite{RSS_WP2}]
  $H_{-1,1}(\riem) \subseteq \Omega_{-1,1}(\riem)$ and inclusion is holomorphic.
  Furthermore
  \[  L^\infty_{-1,1}(\riem) \cap L^2_{\mathrm{hyp}}(\riem) = H_{-1,1}(\riem) \oplus (\mathcal{N} \cap L^2_{\mathrm{hyp}}(\riem)).  \]
 \end{theorem}
 \begin{theorem}[\cite{RSS_WP2}]
  Let $\alpha_t$ be any holomorphic curve in $T_{\mathrm{WP}}(\riem)$ for $t$ in a
  disk centered at $0$, such that $\alpha_0 = [0]$.  There is an open disk $D$ centered at zero so that for each $t \in D$,
  $\alpha_t$ has a representative which is in $L^2_{\mathrm{hyp}}(\riem) \cap L^\infty_{-1,1}(\riem)$,
  which is a holomorphic curve in the Hilbert space $L^2_{\mathrm{hyp}}(\riem)$ and the Banach space $L^\infty_{-1,1}(\riem)$.
 \end{theorem}
 Thus the tangent space at the identity is described by $H_{-1,1}(\riem)$.
 This can be right translated, and also leads to a complex structure as in the
 the classical case.

   The harmonic Beltrami differentials also induce local coordinates.
 For $\mu \in H_{-1,1}(\riem)$ let $f_\mu:\riem \rightarrow \riem_1$ denote a quasiconformal
 solution to the Beltrami equation.
 \begin{theorem}[\cite{RSS_WP2}]
  There is an open neighbourhood of $U$ of $0$ in $H_{-1,1}(\riem)$ such that
  the map $\mu \mapsto [\riem,f_\mu,\riem_1] \in T_{\mathrm{WP}}(\riem)$ obtained by solving
  the Beltrami equation
  is a biholomorphism onto its image.  Furthermore, change of base point is a biholomorphism.
  This describes a system of complex coordinates which endows $T_{\mathrm{WP}}(\riem)$
  with a complex Hilbert manifold structure.
 This complex structure is compatible
  with the complex structure obtained from the fiber structure.
 \end{theorem}
 \begin{corollary}[\cite{RSS_WP2}]  $T_{\mathrm{WP}}(\riem)$ has a finite Weil-Petersson pairing on each
  tangent space.
 \end{corollary}

 \begin{remark}
  In particular, this shows
 that if one had defined the complex structure using harmonic Beltrami differentials in the first place, then $T_{\mathrm{WP}}(\riem)$
 would have a holomorphic fiber structure with fibers (mod $\text{DB}$) biholomorphic
 to $\Oqc_{\mathrm{WP}}(\riem^P_1)$ (and hence locally biholomorphic to $({\Oqc_{\mathrm{WP}}})^{n}$).
 The analogous result for the classical $L^\infty$ case is Theorem  \ref{th:fiber_structure}
 above; even in the $L^\infty$ case it is non-trivial.
 \end{remark}

 The above discussion leads naturally to the following problems.
 \begin{question}  Show that the complex structure on $T_{\mathrm{WP}}(\riem)$ induced by harmonic
  Beltrami differentials $H_{-1,1}(\riem)$ and right composition
  is equivalent to that induced by the Bers embedding into hyperbolically $L^2$ quadratic differentials.\footnotemark
 \end{question}
 \begin{question}  \label{qu:extend_to_infinite_genus_WP}
  Can the holomorphic fiber structure of Theorem \ref{th:fiber_structure} be
  extended to more general surfaces satisfying Lehner's condition,
  in the WP-class/$L^p$ case?  For example,
  to surfaces of infinite genus and/or infinitely many boundary curves?
 \end{question}
 \footnotetext{See Footnote 3}

 Since Takhtajan and Teo obtain a foliation of $T(\disk^*)$ by
 translates of $T_{\mathrm{WP}}(\disk^*)$ we are led to ask:
 \begin{question}  Can $T(\riem)$ be endowed with a complex Hilbert manifold structure,
  such that $T_{\mathrm{WP}}(\riem)$ is the connected component of the identity, and the other
  connected components are right translates of $T_{\mathrm{WP}}(\riem)$?
 \end{question}

\end{subsection}
\begin{subsection}{K\"ahler potential of Weil-Petersson metric}  \label{se:Kahler}
 In this section we give a brief overview of
 some geometric problems associated with the Weil-Petersson metric.

 Ahlfors \cite{Ahlfors_Kahler} showed that for compact Riemann
 surfaces the Weil-Petersson metric is K\"ahler.  In the
 commentary to his collected works he stated that Andr\'e Weil also had a proof but had
 not published it.  Later Ahlfors computed the curvatures of holomorphic sections and the
 Ricci curvature \cite{Ahlfors_curvature}.

 Kirillov and Yuriev \cite{KY2} sketched a generalization of the period mapping in Teichm\"uller theory
 to $\text{Diff}(\mathbb{S}^1)/\text{M\"ob}(\mathbb{S}^1)$.  Some aspects were filled out and extended to the full universal
 Teichm\"uller space by S. Nag \cite{Nag_bulletin}
 and S. Nag and D. Sullivan in \cite{NagSullivan}.    In this formulation, the period map is a map from
 Teichm\"uller space into an infinite-dimensional Siegel disk, which is a set of
 bounded, symmetric operators $Z$ such that $I - Z \overline{Z}$ is positive-definite, analogous to the
 period mapping for compact Riemann surfaces.   Nag and Sullivan \cite{NagSullivan}
 indicated that this period map is holomorphic by proving G\^ateaux holomorphicity.
  The first complete proof of holomorphicity of the Kirillov-Yuriev-Nag-Sullivan
  period mapping was given by Takhtajan and Teo \cite{Takhtajan_Teo_Memoirs}, in both the Weil-Petersson
 and classical setting.\\

 In \cite{HongRajeev} D. K. Hong and S. G. Rajeev showed that the Siegel disk possesses a natural K\"ahler metric, whose K\"ahler potential is
 given, up to a multiplicative constant, by $\log{\text{det} (I-Z \overline{Z})}$.
 Kirillov and Yuriev \cite{KY2}
  and Nag \cite{Nag_bulletin} showed that the pull-back
 of a natural K\"ahler metric on the infinite Siegel disk is the Weil-Petersson metric.
 Thus $\log{\text{det} (I-Z \overline{Z})}$ is also a K\"ahler potential
 for the Weil-Petersson metric.     Hong and Rajeev
 noted that $\text{Diff}(\mathbb{S}^1)/\text{M\"ob}(\mathbb{S}^1)$ was not complete with
 respect to the K\"ahler metric, which indicated that it was not the correct analytic setting
 for the Weil-Petersson metric.  \\

 The completion of $\text{Diff}(\mathbb{S}^1)/\text{M\"ob}(\mathbb{S}^1)$ with respect to the
 Weil-Petersson metric is $T_{\mathrm{WP}}(\disk^*)$, as was demonstrated in \cite{Takhtajan_Teo_Memoirs}.
 In \cite{Takhtajan_Teo_Memoirs} Takhtajan and Teo also proved the striking result that the period mapping $Z$ is in fact the Grunsky
 operator.  In particular, this implies that a constant multiple of
 the Fredholm determinant (\ref{defn:fred determinant}) is a K\"ahler potential for the Weil-Petersson
 metric.  Finally, Takhtajan and Teo \cite{Takhtajan_Teo_Memoirs} and Shen \cite{ShenGrunsky}
 independently showed that the Grunsky operator $Z$ is Hilbert-Schmidt if and only if
 the corresponding Teichm\"uller representative is WP-class.  This is a very satisfying result
 because this is exactly the condition required in order for $Z \overline{Z}$ to be trace-class
 and hence that $\det(I- Z \overline{Z})$ exists.\\

 This leads us to the following natural problems:
 \begin{question}  Is the Weil-Petersson metric on the $L^2$ Teichm\"uller space of bordered surfaces of type $(g,n)$
 K\"ahler?
 More generally, is this true for surfaces satisfying Lehner's condition?
 \end{question}
 \begin{question}  Compute the sectional or Ricci curvatures of the Weil-Petersson
  metric for bordered surfaces of type $(g,n)$, or more generally those
  satisfying Lehner's condition.\footnotemark
 \end{question}
  \footnotetext{After this chapter was submitted in May 2016,
 the paper \cite{Yanagishita_Kahler} of M. Yanagishita appeared, in which it was shown that for surfaces satisfying Lehner's condition, the Weil-Petersson metric is indeed K\"ahler and the sectional and Ricci curvatures are negative.
The convergent Weil-Petersson metric was obtained independently of Radnell, Schippers and Staubach \cite{RSS_WP1,RSS_WP2}.  Yanagishita \cite{Yanagishita_Kahler} also showed that the complex structure from harmonic Beltrami differentials is compatible with the complex structure from the Bers embedding.  When combined with the results of \cite{RSS_WP2}, this apparently shows that these two complex structures are equivalent to that obtained from fibrations over the compact surfaces for surfaces of type $(g,n)$.}

 These problems are closely related to the following.
 \begin{question}  Is there a generalization of the Kirillov-Yuriev-Nag-Sullivan
  period mapping to Teichm\"uller spaces of bordered surfaces of type $(g,n)$, or more general surfaces?
  Can one obtain an analogous K\"ahler potential from this period mapping?
 \end{question}
 Note that generalizations of the Grunsky matrix to higher genus surfaces have been obtained
 by K. Reimer and E. Schippers \cite{Reimer_Schippers}, and for genus zero surfaces with $n$ boundary curves by Radnell, Schippers and Staubach \cite{RadnellSchippersStaubach_Grunsky_genuszero}.
\end{subsection}
\begin{subsection}{Applications of Weil-Petersson class Teichm\"uller theory}
 The Weil-Petersson Teichm\"uller space has recently attracted a great deal of attention,
 in part because of its intrinsic importance to Teichm\"uller theory, and in part
 because of its many applications.  We sketch some of these now.

 In potential theory, the Weil-Petersson
 class domains are precisely those on which the Fredholm determinant of function theory exists,
 and therefore on which it is a viable tool in potential theory.  In Teichm\"uller theory,
 the Weil-Petersson metric has been an important tool in the investigation of the geometry
 of Teichm\"uller space, see e.g. S. Wolpert \cite{Wolpert}. It is now available in vastly greater generality.\\

 There are also various physical applications.  It has been suggested by several authors
  that the universal Teichm\"uller space could serve as a basis
 for a non-perturbative formulation of bosonic string theory.
 See Hong and Rajeev \cite{HongRajeev}, M. G. Bowick and S. G. Rajeev
 \cite{BowickRajeev}, and the (somewhat dated) review of O. Pekonen \cite{Pekonen};
  for a review of connections to conformal field theory see Markina and Vasil'ev \cite{Markina_Vasilev_Virasoro}.
  The $g$-loop scattering amplitude in string theory can be expressed as an
  integral over Teichm\"uller space \cite{HongRajeev}, and thus a non-perturbative
  formulation might be given on the universal Teichm\"uller space, since it contains
  all other Teichm\"uller spaces.  Hong and Rajeev also propose
  that the computation of the scattering amplitudes should involve the exponential
  of the K\"ahler potential discussed in the previous section.
   As observed above, it was known already to Hong and Rajeev
  that $\text{Diff}(\mathbb{S}^1)/\text{M\"ob}(\mathbb{S}^1)$ is not complete with
 respect to the K\"ahler metric (which we now know to be the Weil-Petersson metric).
 The correct analytic setting for Hong and Rajeev's proposal thus appears to be the Weil-Petersson
 class Teichm\"uller space.\\

 The Weil-Petersson Teichm\"uller space also has deep connections with
 two-dimensional conformal field theory as formulated by Segal, Kontsevich, Vafa and others
 (see \cite{Huang} for a review of the literature in the formulation of CFT).  Work of
 Radnell, Schippers and Staubach has established that the Weil-Petersson class rigged
 moduli space is the completion of the analytically rigged moduli space of Friedan/Shenker/Segal/Vafa,
 and is also the largest space on which constructions in conformal field theory can be carried out.
 These include for example sewing properties of the determinant line bundle over the rigged moduli
 space and the existence of local holomorphic sections.   A review of this work can be found in
 \cite{RadnellSchippersStaubach_CFT_review}.\\

 Finally, there are applications to fluid mechanics and infinite dimensional groups of
 diffeomorphisms.  The setting
 for this is a deep insight of Arnol'd, namely that the geodesic equations on infinite-dimensional
 diffeomorphism groups are analogous to the Euler equations of fluid mechanics \cite{Khesin_Wendt}.  Different
 choices of groups and metrics lead to different geodesic equations, which
 in turn are different systems of partial differential equations \cite[Table 4.1]{Khesin_Wendt}.   See E. Grong, I. Markina
 and A. Vasiliev  \cite{GrongMarkinaVasilev} for a survey of choices on $\text{Diff}(\mathbb{S}^1)/\text{M\"ob}(\mathbb{S}^1)$
 and their relation to sub-Riemannian geometry.    In the
 case of the Weil-Petersson metric on $\text{Diff}(\mathbb{S}^1)/\mathbb{S}^1$, the geodesic
 equations are related to the KdV equation.  M. E. Schonbek, A. N. Todorov and J. P. Zubelli
 \cite{Schonbek_Todorov_Zubelli}  were able to obtain long-term solutions to the KdV equation
  using this connection.  F. Gay-Balmaz \cite{Gay} was able to obtain global existence and uniqueness of the
  geodesics,  and applied the Euler-Poincar\'e reduction process to obtain the spatial representation of the
 geodesics. Further important and interesting applications in this direction were given by A. Figalli \cite{Fi},
 F. Gay-Balmaz and  T. S. Ratiu \cite{GR}, and S. Kushnarev \cite{KU}.

 We note that an important technical problem in this direction was solved by Shen \cite{Shen_characterization}.
 Much of the analysis in the fluid mechanical models above has involved the assumption that the
 corresponding quasisymmetries were in $H^{3/2 -\epsilon}(\mathbb{S}^1)$ for $\epsilon>0$, and it was
 an open question whether the quasisymmetries
 in the Weil-Petersson class Teichm\"uller space would be precisely those in $H^{3/2}(\mathbb{S}^1)$.  Theorem \ref{th:Shen_characterization} of
 Shen above gave the correct characterization, and in the same paper
 \cite{Shen_characterization} he also showed that there are WP-class quasisymmetries which are not in
 $H^{3/2}(\mathbb{S}^1)$.
\end{subsection}
\end{section}
\begin{section}{Conclusion}
\begin{subsection}{Concluding remarks}
 In this paper, we have given a number of examples of the general phenomenon of comparison
 moduli spaces in geometric function theory and moduli spaces of Riemann surfaces.  We have
 seen that this concept spontaneously arises in both modern and classical complex analysis.
 We have also attempted to illustrate how
 this notion of moduli space captures many complex analytic phenomena in a
 simple way.

 We would like to conclude with another observation.
 One is struck by the pervasive relevance of classical function theory.  We have
 seen, for example, the unwitting
 re-invention in two-dimensional conformal field theory of the Teichm\"uller
 space of bordered surfaces and conformal welding, and the use of Schiffer's variational technique to construct a complex structure on
 Teichm\"uller space and the rigged moduli space of conformal field theory.  We have
 also seen the Fredholm determinant of classical potential theory - as reformulated by
 Schiffer - emerge as a fundamental
 geometric object on moduli spaces of Riemann surfaces.  These geometric problems in turn
 require the formulation and solution of analytic problems which can only be
 approached with function theory. Other examples spring readily
 to the mind of any mathematician with their ear to the ground.

 We cannot express this in any better way than Ahlfors did \cite{Ahlfors_classical_contemporary}: ``We start out
 from a purely classical problem, and we place it in a much more general modern setting, sometimes
 in a form that would not have been available to a classical mathematician.  When the generalized
 problem is analyzed, it turns out to lead forcefully to a new and evidently significant problem
 in the original purely classical framework.  In other words, we are faced with new evidence of
 the scope and fertility of classical analysis.''
\end{subsection}
\end{section}

\end{document}